%% file: projected_tensor_products.tex
\documentclass[final, 11pt]{elsarticle}
\usepackage[margin=1in]{geometry}

\input{macros}

\usepackage{lineno} 
\usepackage{anyfontsize}
\numberwithin{equation}{section}

\journal{Linear Algebra and its Applications}

\begin{document}

\begin{frontmatter}



\title{Projected Tensor-Tensor Products for Efficient Computation of Optimal Multiway Data Representations}

%

\author[label1]{Katherine Keegan}
\author[label1]{Elizabeth Newman\corref{cor1}} 

 \cortext[cor1]{Corresponding Author}
  \ead{elizabeth.newman@emory.edu}
 \ead[https://math.emory.edu/~enewma5/]{https://math.emory.edu/~enewma5/}
 
\affiliation[label1]{organization={Emory University},
            addressline={400 Dowman Drive}, 
            city={Atlanta},
            postcode={30322}, 
            state={Georgia},
            country={United States of America}}

\begin{abstract}
Tensor decompositions have become essential tools for feature extraction and compression of multiway data. Recent advances in tensor operators have enabled desirable properties of standard matrix algebra to be retained for multilinear factorizations. Behind this matrix-mimetic tensor operation is an invertible matrix whose size depends quadratically on certain dimensions of the data.  
As a result, for large-scale multiway data, the invertible matrix can be computationally demanding to apply and invert and can lead to inefficient tensor representations in terms of construction and storage costs. 
In this work, we propose a new projected tensor-tensor product that relaxes the invertibility restriction to reduce computational overhead and still preserves fundamental linear algebraic properties. 
The transformation behind the projected product is a tall-and-skinny matrix with unitary columns, which depends only linearly on certain dimensions of the data, thereby reducing computational complexity by an order of magnitude. 
We provide extensive theory to prove the matrix mimeticity and the optimality of compressed representations within the projected product framework. 
We further prove that projected-product-based approximations outperform a comparable, non-matrix-mimetic tensor factorization. 
We support the theoretical findings and demonstrate the practical benefits of projected products through numerical experiments on video and hyperspectral imaging data. 
All code for this paper is available at \url{https://github.com/elizabethnewman/projected-products.git}. 

\end{abstract}

%


\end{frontmatter}



\input{01_introduction}

\input{02_background}

\input{03_projected_products}

\input{03_01_special_considerations}

	\input{03_02_equivalent_projected_products}

\input{03_03_algebraic_properties}

\input{03_04_generalized_projected_products}

\input{04_eckart_young}

	\input{04_01_tsvdqII}

	\input{04_02_hosvd}

\input{05_numerical_experiments}

\input{05_01_video}

\input{05_02_hyperspectral}

\input{05_03_hosvd}

\input{06_conclusions}

\section{Acknowledgements}

The work by E. Newman was partially supported by the National Science Foundation
(NSF) under grants [DMS-2309751] and [DE-NA0003525] and the work by K. Keegan was supported by the Department
of Energy Computational Science Graduate Fellowship [DE-SC0023112]. 

Sandia National Laboratories is a multimission laboratory managed and operated by National Technology \& Engineering Solutions of Sandia, LLC, a wholly owned subsidiary of Honeywell International Inc., for the U.S. Department of Energy’s National Nuclear Security Administration under contract DE-NA0003525.
This paper describes objective technical results and analysis. Any subjective views or opinions that might be expressed in the paper do not necessarily represent the views of the U.S. Department of Energy or the United States Government.

This material is based upon work supported by the U.S. Department of Energy, Office of Science, Office of Advanced Scientific Computing Research, Department of Energy Computational Science Graduate Fellowship under Award Number(s) DE-SC0023112.

\appendix

\input{07_projected_products_example}

\input{07_numerical_experiment_data}

\bibliographystyle{plainnat}
\bibliography{main}

\end{document}

%% file: macros.tex
\usepackage{amsmath, amsfonts, amsthm} 
\usepackage{tikz, pgfplots, pgflibraryplotmarks, pgfkeys, ifthen} 
\usetikzlibrary{shapes.arrows, arrows, decorations.markings}
\usetikzlibrary{calc, positioning}
\usepackage{graphicx} 
\usepackage{xcolor} 
\usepackage{multicol} 
\usepackage{hyperref} 
\usepackage{framed}
\usepackage{cleveref}
\crefname{appendix}{}{}
\usepackage{xfrac}
\usepackage{mathabx}
\usepackage{wrapfig}
\pgfplotsset{compat = 1.3}
\usepackage{algorithm, algorithmicx, algpseudocode}

\definecolor{EmoryBlue}{RGB}{1, 33, 105} 
\definecolor{EmoryDarkBlue}{RGB}{12, 35, 64} 
\definecolor{EmoryMediumBlue}{RGB}{0, 51, 160} 
\definecolor{EmoryLightBlue}{RGB}{0, 125, 186} 
\definecolor{EmoryYellow}{RGB}{242, 169, 0} 
\definecolor{EmoryGold}{RGB}{181, 133, 0} 
\definecolor{EmoryMetallicGold}{RGB}{132, 117, 78}

\definecolor{mycolor0}{rgb}{1, 1, 1}%
\definecolor{mycolor1}{rgb}{0.00000,0.44700,0.74100}%
\definecolor{mycolor2}{rgb}{0.85000,0.32500,0.09800}%
\definecolor{mycolor3}{rgb}{0.92900,0.69400,0.12500}%
\definecolor{mycolor4}{rgb}{0.49400,0.18400,0.55600}%
\definecolor{mycolor5}{rgb}{0.46600,0.67400,0.18800}%
\definecolor{mycolor6}{rgb}{0.30100,0.74500,0.93300}%

\definecolor{jet1}{rgb}{0.0000,0.0000,0.6667}
\definecolor{jet2}{rgb}{0.0000,0.0000,1.0000}
\definecolor{jet3}{rgb}{0.0000,0.3333,1.0000}
\definecolor{jet4}{rgb}{0.0000,0.6667,1.0000}
\definecolor{jet5}{rgb}{0.0000,1.0000,1.0000}
\definecolor{jet6}{rgb}{0.3333,1.0000,0.6667}
\definecolor{jet7}{rgb}{0.6667,1.0000,0.3333}
\definecolor{jet8}{rgb}{1.0000,1.0000,0.0000}
\definecolor{jet9}{rgb}{1.0000,0.6667,0.0000}
\definecolor{jet10}{rgb}{1.0000,0.3333,0.0000}
\definecolor{jet11}{rgb}{1.0000,0.0000,0.0000}

\definecolor{mygreen}{RGB}{111, 113, 88}

\usepackage[most]{tcolorbox}
\usepackage{cleveref}

\theoremstyle{definition}
\newtheorem{definition}{Definition}[section]
\newtheorem{myrem}{Remark}[section]

\makeatletter
\crefformat{tcb@cnt@mytheo}{theorem~#2#1#3}
\Crefformat{tcb@cnt@mytheo}{Theorem~#2#1#3}
\crefformat{tcb@cnt@mylemma}{lemma~#2#1#3}
\Crefformat{tcb@cnt@mylemma}{Lemma~#2#1#3}
\crefformat{tcb@cnt@mycorollary}{corollary~#2#1#3}
\Crefformat{tcb@cnt@mycorollary}{Corollary~#2#1#3}
\crefformat{tcb@cnt@myproto}{prototype problem~#2#1#3}
\Crefformat{tcb@cnt@myproto}{Prototype Problem~#2#1#3}
\crefformat{tcb@cnt@myalg}{algorithm~#2#1#3}
\Crefformat{tcb@cnt@myalg}{Algorithm~#2#1#3}
\crefformat{tcb@cnt@myexam}{example~#2#1#3}
\Crefformat{tcb@cnt@myexam}{Example~#2#1#3}
\makeatother

\newtcbtheorem[auto counter, number within=section]{mytheo}{Theorem}%
{colback=EmoryBlue!5,colframe=EmoryBlue, fonttitle=\bfseries}{thm}

\newtcbtheorem[auto counter, number within=section]{mylemma}{Lemma}%
{colback=gray!5,colframe=gray, fonttitle=\bfseries}{lem}

\newtcbtheorem[auto counter, number within=section]{mycorollary}{Corollary}%
{colback=gray!5,colframe=gray, fonttitle=\bfseries}{cor}

\newtcbtheorem[auto counter, number within=section]{myalg}{Algorithm}%
{colback=white!5,colframe=black, fonttitle=\bfseries}{alg}

\newtcbtheorem[auto counter, number within=section]{myexam}{Example}%
{colback=white!5,colframe=EmoryGold, fonttitle=\bfseries}{exam}

\usepackage{xcolor}
\usepackage{amsthm}
\usepackage{framed}

\theoremstyle{plain}

\theoremstyle{definition}
\newtheorem{protoactionitems}{Action Items}[section]

\newcommand{\starM}{\star_{\bfM}}
\newcommand{\starQ}{\star_{\bfQ^H}^\prime}
\newcommand{\starW}{\star_{\bfW^H}^\prime}
\newcommand{\starQperp}{\star_{\bfQ_{\perp}^H}^\prime}

\DeclareMathOperator*{\rank}{rank}

\DeclareMathOperator*{\argmax}{arg\ max}
\DeclareMathOperator*{\argmin}{arg\ min}
\DeclareMathOperator{\trank}{\starM-rank}

\DeclareMathOperator{\projtrank}{\starQ-rank}

\DeclareMathOperator*{\myVec}{vec}
\DeclareMathOperator{\myFold}{fold}

\DeclareMathOperator{\Stiefel}{St}
\DeclareMathOperator{\GL}{GL}

\newcommand{\blue}[1]{{\color{blue} #1}}

\definecolor{mycolor1}{rgb}{0.00000,0.44700,0.74100}%
\definecolor{mycolor2}{rgb}{0.85000,0.32500,0.09800}%
\definecolor{mycolor3}{rgb}{0.92900,0.69400,0.12500}%
\definecolor{mycolor4}{rgb}{0.49400,0.18400,0.55600}%
\definecolor{mycolor5}{rgb}{0.46600,0.67400,0.18800}%
\definecolor{mycolor6}{rgb}{0.30100,0.74500,0.93300}%

\def\mydefb#1{\expandafter\def\csname bf#1\endcsname{\mathbf{#1}}}
\def\mydefallb#1{\ifx#1\mydefallb\else\mydefb#1\expandafter\mydefallb\fi}
\mydefallb aAbBcCdDeEfFgGhHiIjJkKlLmMnNoOpPqQrRsStTuUvVwWxXyYzZ\mydefallb

\def\mydefb#1{\expandafter\def\csname #1bb\endcsname{\mathbb{#1}}}
\def\mydefallb#1{\ifx#1\mydefallb\else\mydefb#1\expandafter\mydefallb\fi}
\mydefallb aAbBcCdDeEfFgGhHiIjJkKlLmMnNoOpPqQrRsStTuUvVwWxXyYzZ\mydefallb

\def\mydefb#1{\expandafter\def\csname #1cal\endcsname{\mathcal{#1}}}
\def\mydefallb#1{\ifx#1\mydefallb\else\mydefb#1\expandafter\mydefallb\fi}
\mydefallb aAbBcCdDeEfFgGhHiIjJkKlLmMnNoOpPqQrRsStTuUvVwWxXyYzZ\mydefallb

\def\mydefb#1{\expandafter\def\csname T#1\endcsname{\boldsymbol{\mathcal{#1}}}}
\def\mydefallb#1{\ifx#1\mydefallb\else\mydefb#1\expandafter\mydefallb\fi}
\mydefallb aAbBcCdDeEfFgGhHiIjJkKlLmMnNoOpPqQrRsStTuUvVwWxXyYzZ\mydefallb

\def\mydefgreek#1{\expandafter\def\csname bf#1\endcsname{\text{\boldmath$\mathbf{\csname #1\endcsname}$}}}
\def\mydefallgreek#1{\ifx\mydefallgreek#1\else\mydefgreek{#1}%
   \lowercase{\mydefgreek{#1}}\expandafter\mydefallgreek\fi}
\mydefallgreek {alpha}{Alpha}{beta}{Beta}{gamma}{Gamma}{delta}{Delta}{epsilon}{Epsilon}{zeta}{Zeta}{eta}{Eta}{theta}{Theta}{iota}{Iota}{kappa}{Kappa}{lambda}{Lambda}{mu}{Mu}{nu}{Nu}{omicron}{Omicron}{pi}{Pi}{rho}{Rho}{sigma}{Sigma}{tau}{Tau}{upsilon}{Upsilon}{phi}{Phi}{xi}{Xi}{chi}{Chi}{psi}{Psi}{omega}{Omega}\mydefallgreek

\usepackage{subcaption}

\definecolor{purple1}{RGB}{135, 225, 0} 
\definecolor{purple2}{RGB}{193, 252, 129} 
\definecolor{purple3}{RGB}{228, 244, 210} 

\definecolor{teal1}{RGB}{80, 227, 194} 
\definecolor{teal2}{RGB}{167, 246, 228} 
\definecolor{teal3}{RGB}{221, 247, 239}

\definecolor{green1}{RGB}{234, 140, 229} 
\definecolor{green2}{RGB}{235, 177, 237} 
\definecolor{green3}{RGB}{246, 222, 244}

%% file: 01_introduction.tex
\section{Introduction}
\label{sec:introduction}
Multiway arrays or tensors arise naturally across modern data science applications, such as precision medicine~\cite{Mor2022, LarssonEtAl2007}, signal processing~\cite{7891546}, and machine learning~\cite{NewmanEtAl2024, 10.5555/2969239.2969289}. Tensor decompositions, typically framed as high-dimensional analogs of the matrix singular value decomposition (SVD), have become widely used to efficiently represent multiway data for subsequent computation and analysis~\cite{TensorTextbook, KoldaBader}.  Tensor factorizations come in many varieties, from the classical Canonical Polyadic/Parallel Factor (CP) decomposition~\cite{Hitchcock1927, Harshman1970, CarrollChang1970} and Tucker decomposition~\cite{Tucker1966, LathauwerMoorVandewalle2000} to the more recent tensor train and tensor network decompositions~\cite{Oseledets2011, Cichocki_2016}.  Modern advancements of multilinear decompositions exploit underlying structure cleverly to provide theoretical insights and accelerate computation, including using tools from algebraic geometry to decompose symmetric tensors~\cite{WANG202369, kileel2024subspacepowermethodsymmetric} and incorporating randomized sketching for efficient implementation and storage of high-dimensional data~\cite{NEURIPS2022_fe91414c, NEURIPS2018_45a766fa}.  

Our paper focuses on building new computational advancements while retaining algebraic advantages from the $\starM$-framework (the prefix is pronounced ``star-M'' or ``M''), a matrix-mimetic framework that views tensors as operators~\cite{KernfeldKilmer2015}. The multilinear operation, called the $\starM$-product, multiplies two tensors under an algebraic ring operation determined by an invertible matrix $\bfM$. As a result, the $\starM$-product ``looks and feels'' like matrix multiplication, and thereby preserves familiar linear algebra properties. In particular, a tensor SVD under the $\starM$-product yields provably optimal compressed representations that can theoretically and empirically outperform the matrix SVD and comparable tensor factorizations~\cite{kilmer2019tensortensor}. The optimality of the representations is the hallmark of the $\starM$-framework; other classical tensor decompositions only achieve quasi-optimality. Recent work has leveraged $\starM$-optimality to optimize the choice of invertible $\bfM$ and further improve the quality of the $\starM$-representations~\cite{newman2024optimalmatrixmimetictensoralgebras}.  

For sufficiently large multiway data, a computational bottleneck of the $\starM$-product is the storage and application of $\bfM$ and its inverse. Remedies include using easy-to-invert structure of $\bfM$ (e.g., unitary) and storing the matrix implicitly by, e.g., using the fast (inverse) Fourier transform, as in the original $t$-product~\cite{KilmerMartin2011, KilmerBramanHaoHoover2013}.  
Even with these remedies, the requirement of invertibility prevents $\starM$-representations from compressing along certain dimensions or modes of the data, which can lead to prohibitively expensive computational and storage costs.

\subsection{Our Contributions}
\label{sec:our_contributions}

In this work, we introduce a new projected tensor-tensor product as a practical relaxation of the original $\starM$-product. 
Our new tensor-tensor product is defined by a matrix $\bfQ$ with unitary columns that is not necessarily invertible. 
This choice of $\bfQ$ reduces the computational complexity and representation storage costs by an order of magnitude based on the size of the multiway data. Our contributions include developing new and extensive theoretical foundations for projected products, including proofs that the projected product is matrix mimetic. Notably, we achieve Eckart-Young optimality results under the projected product algebra and provide insight into an optimal choice of projected product matrix $\bfQ$. We further prove that representations under the projected product yield better approximations than the higher-order SVD (HOSVD)~\cite{LathauwerMoorVandewalle2000}. Our numerical experiments provide strong empirical support of our theoretical results and demonstrate the ability of the projected product representations to approximate multiway data well with significant storage reduction. For transparency and reproducibility, we provide our an open-source at \url{https://github.com/elizabethnewman/projected-products.git}.

\subsection{Related Work}
\label{sec:related_work}

Other non-invertible tensor-tensor products been proposed for the specific application of tensor completion via tensor nuclear norm minimization~\cite{Zhang_2014_CVPR}. In~\cite{9115254}, the authors introduce a framelet transform as an alternative to the fast Fourier transform (\texttt{fft}). aSimilar to the \texttt{fft}, framelet transforms can be implemented implicitly and efficiently with little additional storage overhead. However, the proposed framelet transformation increases the dimensions of the tensor during application, resulting in greater computational cost within the tensor completion algorithm. Similarly, in~\cite{9525838}, the authors propose a dictionary-based transformation to ideally produce sparse tensor representations. However, the dictionary must be stored as an overcomplete matrix, which ultimately increases the computational and storage costs if the solution is insufficiently sparse. In~\cite{KongLuLin2021}, the authors propose learning a data-dependent transformation with orthonormal columns (semi-orthogonal) as a subproblem of tensor nuclear norm minimization.  The paper focuses on developing two strategies to optimize the matrix, variance maximization and manifold optimization, without developing nor leveraging algebraic properties that the underlying product induces. The work in~\cite{10.1007/s10915-022-01937-1} extends from~\cite{KongLuLin2021} by applying a semi-orthogonal matrix followed by pointwise nonlinearity and learning the matrix through an alternating minimization strategy.  The nonlinearity is generalized further in~\cite{9780890}, which trains a multi-layer neural network as the transformation. This design introduces new flexibility to the tensor-tensor product, but also potentially increases storage, computation, and training costs.  

While the works of the above papers present practical advancements of non-invertible tensor-tensor products for tensor completion, none provide insight into the the algebraic implications of an underlying non-invertible transformation.  Our work develops a unified algebraic framework for tensor algebras defined by real- or complex-valued matrices with unitary columns, introduces new theory about the optimality of the tensor representations under these non-invertible products, and extends the types of applications to which this projected product can be applied. 

\subsection{Outline of the Paper}
\label{sec:outline}

The paper proceeds as follows. In~\Cref{sec:background}, we describe the notation and algebraic foundations of the $\starM$-product. In~\Cref{sec:projected_products}, we introduce the projected tensor-tensor product, describe its differences from and relationships to the original $\starM$-product (\Cref{sec:special_considerations} and~\Cref{sec:equivalent_presentations}), and verify its matrix mimeticity (\Cref{sec:algebraic_properties}). In~\Cref{sec:tsvdq}, we present the projected-product-based tensor SVD and a compressible variant (\Cref{sec:tsvdqII}), prove the Eckart-Young optimality for both representations, and verify that the tensor SVD representations can outperform the higher-order SVD (\Cref{sec:hosvd}). In~\Cref{sec:numericalexperiments}, we empirically support the theoretical results and demonstrate the high-quality, compressed representations we can obtain using projected products through several numerical experiments on both video and hyperspectral imaging data. 
We conclude in~\Cref{sec:conclusions} with a discussion of future directions.

%% file: 02_background.tex
\section{Background}
\label{sec:background}

Tensors, denoted in bold calligraphic letters $\TA$, are multiway arrays and the order of a tensor is the number of dimensions or modes. Familiar linear algebra objects can be interpreted as tensors; scalars, denoted with lowercase $a \in \Cbb$, are order-$0$ tensors; vectors, denoted with bold lowercase $\bfa \in \Cbb^{n_1}$, are order-$1$ tensors; and matrices, denoted with bold uppercase $\bfA \in \Cbb^{n_1\times n_2}$, are order-$2$ tensors. We use the term \emph{tensor} to refer to arrays of order-$3$ or higher.  This paper will focus on order-$3$ tensors, though higher-order extensions can be made recursively; see~\cite{keegan2022a}. 

\begin{figure}
\centering

\subfloat[\begin{tabular}{c} Tensor \\ $\TA \in \Cbb^{n_1\times n_2\times n_3}$ \end{tabular}]{\includegraphics[width=0.225\linewidth]{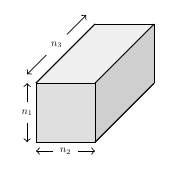}}
\subfloat[\begin{tabular}{c} Lateral slices \\ $\TA_{:,j,:} \in \Cbb^{n_1\times 1\times n_3}$ \end{tabular}]{\includegraphics[width=0.225\linewidth]{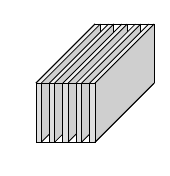}}
\subfloat[\begin{tabular}{c} Frontal slices \\ $\TA_{:,:,k} \in \Cbb^{n_1\times n_2}$ \end{tabular}]{\includegraphics[width=0.225\linewidth]{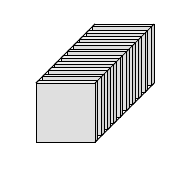}}
\subfloat[\begin{tabular}{c} Tubes\\ $\TA_{i,j,:}\in \Cbb^{1\times 1\times n_3}$ \end{tabular}]{\includegraphics[width=0.225\linewidth]{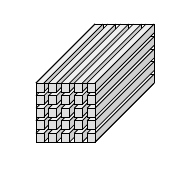}}

\caption{Illustration of key partitions of third-order tensors.}
\label{fig:tensor_notation}
\end{figure}

Analogous to rows and columns of matrices, a tensor can be indexed in various ways. We will use {\sc Matlab} indexing notation to discuss key partitions; e.g., $\bfA_{:,j}$ or $\bfA(:,j)$ indicates the $j$-th column of a matrix. Two key tensor partitions are slices (one index fixed) and fibers (two indices fixed). Let $\TA\in \Cbb^{n_1\times n_2\times n_3}$ be an order-$3$ tensor. Frontal slices $\TA_{:,:,k} \in \Cbb^{n_1\times n_2}$ for $k=1,\dots, n_3$ are matrices stacked along the third dimension and lateral slices $\TA_{:,j,:} \in \Cbb^{n_1\times 1\times n_3}$ for $j=1,\dots,n_2$ are matrices oriented along the third dimension and stacked along the second dimension. Tube fibers $\TA_{i,j,:}\in \Cbb^{1\times 1\times n_3}$ for $i=1,\dots,n_1$ and $j=1,\dots,n_2$ are vectors the lying along the third dimension. We provide an illustration of the various tensor partitions in~\Cref{fig:tensor_notation}. The Frobenius norm of an order-$3$ tensor can be defined via the frontal slices; that is, $\|\TA\|_F^2 =  \sum_{k=1}^{n_3}\|\TA_{:,:,k}\|_F^2$.

To operate on tensors, we define the mode-$3$ and facewise products.  Extensions to order-$d$ tensors can be found in~\cite{KoldaBader}. 

\begin{definition}[mode-$3$ unfolding/folding]\label{def:mode3unfolding}
The mode-$3$ unfolding of $\TA \in \Cbb^{n_{1} \times n_{2} \times n_{3}}$, denoted $\bfA_{(3)}\in \Cbb^{n_3\times n_1n_2}$, matricizes the tensor such that the columns are vectorized tubes. Specifically, 
    	\begin{align}
	\bfA_{(3)}(:,J(i,j)) = \myVec(\TA_{i,j,:}) \qquad \text{where} \qquad J(i,j) = i + (j-1)n_1
	\end{align}
for $i=1,\dots,n_1$ and $j=1,\dots,n_2$. Here, $\myVec: \Cbb^{1\times 1\times n_3} \to \Cbb^{n_3}$ turns a tube into a column vector. The mode-$3$ folding, $\myFold_{(3)}$,  reverses mode-$3$ matricization. 
\end{definition}

\begin{definition}[mode-$3$ product]
The mode-$3$ product between $\TA \in \Cbb^{n_{1} \times n_{2} \times n_{3}}$ and $\bfM \in \Cbb^{p \times n_{k}}$, denoted $\TA \times_3 \bfM\in \Cbb^{n_1\times n_2\times p}$, is given by 
    \begin{align}
    \TA \times_{3} \bfM =\myFold_{(3)}(\bfM \bfA_{(3)}).
    \end{align}
\end{definition}

\begin{definition}[facewise product]
The facewise product between $\TA \in \Cbb^{n_{1} \times m \times n_{3}}$ and $\TB \in \Cbb^{m \times n_{2} \times n_{3}}$, denoted $\TA \smalltriangleup \TB \in \Cbb^{n_1 \times n_2\times n_3}$, multiplies the corresponding frontal slices together; i.e.,
    \begin{align}
    (\TA \smalltriangleup \TB)_{:,:,i} = \TA_{:,:,i}\TB_{:,:,i}  \qquad \text{for $i=1,...,n_{3}$.}
    \end{align}
\end{definition}

Combining the mode-$3$ and facewise products, we now define the $\starM$-product as our foundational tensor-tensor product. As a shorthand, we will write that a complex-valued, invertible $n_3\times n_3$ matrix $\bfM$ belongs to the general linear group; that is, $\bfM \in \GL_{n_3}(\Cbb)$. 

\begin{definition}[$\starM$-product]
Let $\TA \in \Cbb^{n_{1} \times m \times n_{3}}$ and $\TB \in \Cbb^{m \times n_{2} \times n_{3}}$ and let $\bfM \in \GL_{n_3}(\Cbb)$. The $\starM$-product is given by 
    \begin{align}
    \TA \starM \TB = \left ( (\TA \times_{3} \bfM)\smalltriangleup (\TB \times_{3} \bfM ) \right )\times_{3} \bfM^{-1}.
    \end{align} 
\end{definition} 

We say $\TA$ lies in the \emph{spatial} or data domain. When we apply the transformation $\bfM$ along the tubes, we say $\widehat{\TA} = \TA \times_{3} \bfM$ lies in \emph{transform} or frequency domain. We denote tensors in the transform domain with the ``hat'' notation. The origins of this terminology and notation come from the original $t$-product~\cite{KilmerMartin2011}, which used  the (unnormalized) discrete Fourier transform as the transformation matrix $\bfM$. Hence, we adopt the term ``frequency'' or ``transform'' domain for general  transformation matrices. 

The cornerstone of the $\starM$-framework is its matrix mimeticity, which naturally extends properties from standard matrix multiplication to tensors. We see evidence of matrix mimeticity a tube-wise presentation of the $\starM$-product
	\begin{align}\label{eq:mprodEntrywise}
	(\TA \starM \TB)_{i,j,:} = \sum_{\ell=1}^m \TA_{i,\ell,:} \starM \TB_{\ell,j,:}
	\end{align}
for $i=1,\dots,n_1$ and $j=1,\dots,n_2$. This is analogous to the entrywise definition of matrix-matrix multiplication where tubes act as scalars.

%% file: 03_projected_products.tex
\section{Projected Tensor-Tensor Products}
\label{sec:projected_products}

A major restriction of the $\starM$-product is that the transformation matrix $\bfM$ has to be invertible. While this yields algebraic advantages, the computational and storage costs of resulting representations can, in some cases,  be dominated by $\bfM$. We introduce a new projected tensor-tensor product that significantly reduces the computational overhead of applying and storing the transformation. Our key modification is to consider transformation matrices $\bfQ \in \Cbb^{n_3\times p}$ that have unitary columns, but are not necessarily invertible. In the language of manifolds, we say $\bfQ$ belongs to the Stiefel manifold over the complex numbers; that is, $\bfQ \in \Stiefel_{n_3,p}(\Cbb)$ where $\Stiefel_{n_3,p}(\Cbb) = \{\bfX \in \Cbb^{n_3\times p} \mid \bfX^H \bfX = \bfI_p\}$.

\begin{figure}
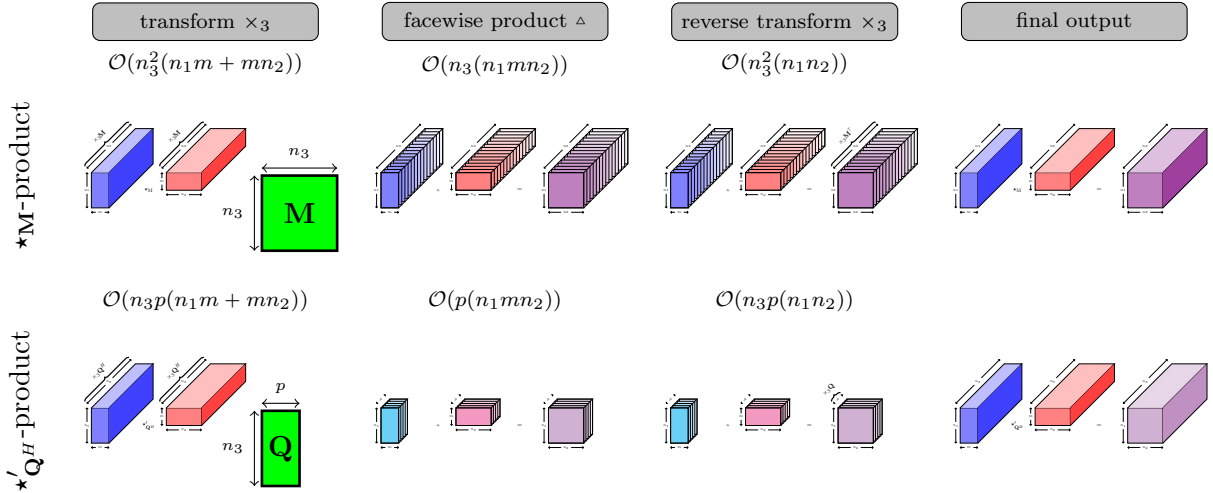


\centering

\begin{tikzpicture}

\def\w{3.75}
\def\s{0.1}

\def\hh{\vphantom{$\Ocal(n_3^2 (n_1m + mn_2))$}}
\foreach[count=\a] \i/\header/\cost/\projcost in {
	2/{transform $\times_3$}/$\Ocal(n_3^2 (n_1m + mn_2))$/$\Ocal(n_3p (n_1m + mn_2))$, 
	4/{facewise product $\smalltriangleup$}/$\Ocal(n_3 (n_1m n_2))$/$\Ocal(p (n_1m n_2))$, 
	5/{reverse transform $\times_3$}/$\Ocal(n_3^2(n_1n_2))$/$\Ocal(n_3p(n_1n_2))$,
	6/{final output}/\phantom{}/~}
	{
	
	\node at (\a*\w+\a*\s,0) (n\a) {\includegraphics[width=\w cm]{projected_products_full\i.pdf}};

	\node[above=0.0cm of n\a.north, anchor=south] (oo) {\scriptsize \cost};
	\node[above=0.6cm of n\a.north, anchor=south, fill=lightgray, draw, rounded corners=0.1cm, minimum width=3cm, minimum height=0.5cm] {\scriptsize \header};
	
	\node[below=0.75cm of n\a.south, anchor=north] (m\a) {\includegraphics[width=\w cm]{projected_products_small\i.pdf}};
	
	\node[above=0.0cm of m\a.north, anchor=south] (oo) {\scriptsize \projcost};

}

\draw[line width=1pt, fill=green] ($(n1) + (1.25+0.5,-0.5-0.5)$) rectangle ($(n1) + (1.25-0.5,-0.5+0.5)$) node[midway] {$\bfM$};
\draw[<->]  ($(n1) + (1.25-0.6,-0.5+0.5)$) -- node[left, midway] {\tiny $n_3$}  ($(n1) + (1.25-0.6,-0.5-0.5)$);
\draw[<->] ($(n1) + (1.25+0.5,-0.5+0.6)$) -- node[above, midway] {\tiny $n_3$} ($(n1) + (1.25-0.5,-0.5+0.6)$);

\draw[line width=1pt, fill=green] ($(m1) + (0.75+0.5,-0.5-0.5)$) rectangle ($(m1) + (1.25-0.5,-0.5+0.5)$) node[midway] {$\bfQ$};
\draw[<->]  ($(m1) + (1.25-0.6,-0.5+0.5)$) -- node[left, midway] {\tiny $n_3$}  ($(m1) + (1.25-0.6,-0.5-0.5)$);
\draw[<->] ($(m1) + (0.75+0.5,-0.5+0.6)$) -- node[above, midway] {\tiny $p$} ($(m1) + (1.25-0.5,-0.5+0.6)$);

\node[rotate=90, above=0.0cm of n1.west, anchor=south] {$\starM$-product\vphantom{$\starQ$}};
\node[rotate=90, above=0.0cm of m1.west, anchor=south] {$\starQ$-product};

\end{tikzpicture}
\caption{Comparison of $\starM$-pipeline (top) and $\starQ$-pipeline (bottom) for multiplying tensors. 
Above each operation, we describe the computational cost for dense numerical linear algebra operations with an easy-to-invert matrix $\bfM$ (see~\cite[Section 1.4.1]{Halko2011}). In this setting, the cost of the $\starM$-product depends quadratically on $n_3$, i.e., $\Ocal(n_3^2)$, whereas the $\starQ$-product has only a linear dependence, i.e., $\Ocal(n_3)$. We note that inverting a general matrix $\bfM$ could increase the $\starM$-product reverse transform cost by a factor of $n_3$. Conversely, if $\bfM$ could be implemented via a fast transformation (e.g., \texttt{fft}), then the cost of the of transforms could decrease to $\Ocal(n_3\log n_3)$.}
\label{fig:starM_vs_starQ}
\end{figure}

\begin{definition}[$\starQ$-product]\label{def:projprod}
Let $\TA \in \Cbb^{n_{1} \times m \times n_{3}}$ and $\TB \in \Cbb^{m \times n_{2} \times n_{3}}$ and let $\bfQ \in \Stiefel_{n_3,p}(\Cbb)$.  Then, the projected tensor-tensor product is defined as 
	\begin{align}
	\TA \starQ \TB = [ (\TA \times_{3} \bfQ^H )\smalltriangleup (\TB \times_{3} \bfQ^H) ]\times_{3} \bfQ.
	\end{align}
\end{definition} 
For ease of discussion, we will equivalently call this the projected product or $\starQ$-product (where the prefix is pronounced ``star-Q prime,'' ``star-Q,'' or ``Q''). We will slightly abuse notation and use the ``hat'' notation to denote a tensor in the transform domain for the projected product; i.e., $\widehat{\TA} = \TA \times_{3} \bfQ^{H}$. Whether the ``hat'' refers to the $\starM$- or $\starQ$-transformation will be clear from context or explicitly stated. We depict the differences between computational and storage costs of $\starM$-product and $\starQ$-product in~\Cref{fig:starM_vs_starQ}.  The main takeaway is that the cost of the $\starM$-product depends quadratically on $n_3$ whereas the $\starQ$-product only depends linearly on $n_3$.

 We call~\Cref{def:projprod} a ``projected product'' because applying the transformation and returning to the spatial domain results in an orthogonal projection\footnote{Here, although the matrices are complex-valued, we use the more common term ``orthogonal projection.''} of the tensor tubes onto the column space of $\bfQ$.  For example, if we apply the transform, $\bfQ^H$, and its pseudoinverse, $\bfQ$, to a tensor $\TA$, we obtain $(\TA \times_3 \bfQ^H) \times_3 \bfQ  = \TA \times_3 \bfQ \bfQ^H \not= \TA$.

We define $\starQ$-versions of linear algebraic concepts, including identity, transposition, unitary, and diagonal. 

\begin{definition}[$\starQ$-identity tensor]\label{def:identity_tensor}
Given $\bfQ \in \Stiefel_{n_3,p}(\Cbb)$, the identity tensor $\TI \in \mathcal{R}^{m \times m \times n_{3}}$ is constructed such that each frontal slice in the transform domain is the identity matrix; that is, 
	\begin{align}
	\widehat{\TI}_{:,:,i} = \bfI_{m} \qquad \text{for $i=1,\dots,p$.}
	\end{align}
\end{definition}

\begin{definition}[$\starQ$-conjugate transpose]\label{def:conjugate_transpose}
Given $\bfQ \in \Stiefel_{n_3,p}(\Cbb)$ and  $\TA\in \Cbb^{n_1\times n_2\times n_3}$, its conjugate or Hermitian transpose $\TA^H \in\Cbb^{n_2\times n_1\times n_3}$ is formed by computing the matrix conjugate transpose of each frontal slice in the transform domain; that is, 
	\begin{align}
	\widehat{(\TA^H)}_{:,:,i} = (\widehat{\TA}_{:,:,i})^H \qquad \text{for $i=1,\dots,p$.}
	\end{align}
\end{definition}

\begin{definition}[$\starQ$-unitary]
Given $\bfQ\in \Stiefel_{n_3,p}(\Cbb)$, we say $\TU \in \Cbb^{m\times m \times n_{3}}$ is unitary if 
	\begin{align}
	\TU^H \starQ \TU = \TU \starQ \TU^H = \TI. 
	\end{align}
\end{definition}

\begin{definition}[facewise diagonal (f-diagonal)] 
A tensor $\TD\in \Cbb^{n_1\times n_2\times n_3}$ is facewise diagonal if its only nonzero entries are contained within its diagonal tubes; that is, $\TD_{i,i,:}$ is potentially nonzero for $i=1,\dots,\min(n_1,n_2)$ and the remaining tubes are zero.
\end{definition}

%% file: 03_01_special_considerations.tex
\subsection{Special Considerations for Projected Products}
\label{sec:special_considerations}

Relaxing the invertibility restriction does have some notable consequences for the uniqueness of algebraic properties. For example, consider the $\starQ$-identity tensor $\TI = \widehat{\TI} \times_3 \bfQ$ where $\widehat{\TI}_{:,:,i} = \bfI_m$ for $i=1,\dots,p$.  Because $\bfQ$ is not invertible, there are infinitely many tensors that are equivalent to the $\starQ$-identity tensor. Specifically, any tensor of the form
	\begin{align}\label{eq:identityNull}
	\TJ = \TI + \TE \times_3 (\bfI_{n_3} - \bfQ\bfQ^H)
	\end{align}
will be an identity tensor.  The tubes of the second term  lie in the null space of $\bfQ^H$, and hence become zero in the transform domain; that is,  $[\TE \times_3 (\bfI_{n_3} - \bfQ\bfQ^H)]\times_3 \bfQ^H  = \bf0$. 

A similar lack of uniqueness can be found for the $\starQ$-conjugate transpose.  For example, if $\TB = \TA + \TE \times_3  (\bfI_{n_3} - \bfQ\bfQ^H)$, then $\TB^H = \TA^H$, but $(\TB^H)^H \not= \TB$. Thus, the $\starQ$-conjugate transpose is not injective. While such nuances of the $\starQ$-product sacrifice some uniqueness properties of the $\starM$-product, the core algebraic concepts are preserved. 

We note another subtle difference between the $\starM$- and $\starQ$-products for the $t$-product specifically. The $t$-product uses $\bfM = \bfF$, the discrete Fourier transform, which consists of entries based on the complex roots of unity. Conveniently, if $\TA$ and $\TB$ are real-valued tensors, $\TA \starM \TB$ will also be real-valued. However, under the projected product with $\bfQ = \bfF_{1:p,:}^H$, $\TA \starQ \TB$ could be complex-valued.  While this is not problematic theoretically, this is an important consideration in practice, particularly when considering storage costs of complex numbers.

%% file: 03_02_equivalent_projected_products.tex
\subsection{Equivalent Presentations of Projected Products}
\label{sec:equivalent_presentations}

There are several equivalent ways to connect projected products with the original $\starM$-product. First, if $\bfM = \bfQ^H$, then the $\starM$- and $\starQ$-products are equal. More generally, let $\bfM \in \Stiefel_{n_3,n_3}(\Cbb)$ and let $\bfQ = \bfM_{1:p,:}^H$ and $\bfQ_{\perp} = \bfM_{p+1:n_3,:}^H$; alternatively, we write $\bfM = \begin{bmatrix} \bfQ & \bfQ_{\perp} \end{bmatrix}^H$.  Note that we can always reorder the rows of $\bfM$ to obtain this partition because of the $\starM$-product invariance to row permutations~\cite[Theorem 2.2]{newman2024optimalmatrixmimetictensoralgebras}.

We first connect the $\starM$-product with the $\starQ$-product through projections via
	\begin{align}\label{eq:colSpaceDef}
	\TA \starQ \TB \equiv (\TA \times_3 \bfQ\bfQ^H) \starM (\TB \times_3 \bfQ\bfQ^H). 
	\end{align}
In essence, the $\starQ$-product is equal to the $\starM$-product after orthogonally-projecting the tensor tubes onto the column space of $\bfQ$. A similar observation can be made for the remaining columns of $\bfM$ using the projection $\bfQ_{\perp} \bfQ_{\perp}^H = \bfI_{n_3} - \bfQ\bfQ^H$. 

A related perspective comes from recognizing that any tensor can be decomposed as the sum 
	\begin{align}\label{eq:tensorSplitProjection}
	\TC = \TC \times_3 \bfQ\bfQ^H + \TC \times_3 (\bfI_{n_3} - \bfQ\bfQ^H).
	\end{align}
The $\starM$-product can thus be expressed as a projected product plus an error term; that is,
	\begin{align}\label{eq:nullspacePerspective}
	\TA \starM \TB 
	&=\TA \starQ \TB + (\TA \starM \TB) \times_3 (\bfI_{n_3} - \bfQ \bfQ^H).
	\end{align}
The tubes of the second term lie in the null space of $\bfQ^H$. Alternatively, because the second term lies in the column space of $\bfQ_{\perp}$, we can express the $\starM$-product as the sum of two projected products
	\begin{align}\label{eq:sum_proj_prod}
	\TA \starM \TB &= \TA \starQ \TB + \TA \starQperp \TB.
	\end{align}
There is a subtlety to this definition. In the spatial domain, the two tensors are the same size and summable. The $\starQ$- and $\starQperp$-transformed tensors may have a different number of frontal slices.  Collectively, those slices form all frontal slices in the $\starM$-transform domain; that is,
	\begin{subequations}\label{eq:proj_prod_split_frontal_slices}
	\begin{alignat}{4}
	\widehat{(\TA \starM \TB)}_{:,:,1:p} &= (\TA \starM \TB) \times_3 \bfM_{1:p,:} &&= (\TA \starQ \TB) \times_3 \bfQ^H &&= \widehat{(\TA \starQ \TB)} \label{eq:starM_vs_starQ_first} \\
	\widehat{(\TA \starM \TB)}_{:,:,p+1:n_3} &= (\TA \starM \TB) \times_3 \bfM_{p+1:n_3,:} &&= (\TA \starQperp \TB) \times_3 \bfQ_{\perp}^H &&= \widehat{(\TA \starQperp \TB)} 
	\end{alignat} 
	\end{subequations}
Consequently, frontal slices of the $\starQ$- and $\starQperp$-products are zeroed out in the $\starM$-transform domain; i.e.,
	\begin{align}
	(\TA \starQ \TB) \times_3 \bfM_{p+1:n_3,:} ={\bf0} \qquad \text{and} \qquad  
	(\TA \starQperp \TB) \times_3 \bfM_{1:n_3,:} ={\bf0}.  \label{eq:backFrontalSlicesZero}
	\end{align}
For concreteness, we provide an example of the various projected product perspectives in~\Cref{sec:projected_products_example}.

%% file: 03_03_algebraic_properties.tex
\subsection{Algebraic Properties of Projected Products}
\label{sec:algebraic_properties}

We present the core algebraic properties of projected products for tubal multiplication: commutativity, associativity, and distributivity.  We also prove that transposition over the projected product of tensors follows the matrix definitions. 

\begin{mytheo}{Commutivity of the Projected Product for Tubes}{symmetry}
Given $\bfQ\in \Stiefel_{n_3,p}(\Cbb)$, the projected product of tubes is commutative; that is, for any tubes $\bfa, \bfb\in \Cbb^{1\times 1\times n_3}$, we have $\bfa \starQ \bfb = \bfb \starQ \bfa.$
\end{mytheo}

\begin{proof}
By~\Cref{def:projprod}, we have
	\begin{align}
	\bfa \starQ \bfb &= [\widehat{\bfa} \odot \widehat{\bfb}] \times_3\bfQ
	=  [\widehat{\bfb} \odot \widehat{\bfa}] \times_3\bfQ = \bfb \starQ \bfa
	\end{align}
where $\odot$ is the Hadamard pointwise product, which itself is commutative.  
\end{proof}

\begin{mytheo}{Algebraic Properties of the Projected Product}{associatvitity}
Given $\bfQ\in \Stiefel_{n_3,p}(\Cbb)$ and arbitrary tubes $\bfa, \bfb, \bfc\in \Cbb^{1\times 1\times n_3}$, the projected product is
\begin{enumerate}
\item associative, i.e., $(\bfa \starQ \bfb) \starQ \bfc = \bfa \starQ (\bfb \starQ \bfc)$ and

\item distributive over addition, i.e., $(\bfa + \bfb) \starQ \bfc = \bfa \starQ \bfc +  \bfb \starQ \bfc$.

\end{enumerate}

\end{mytheo}

\begin{proof} We complete the proof using~\Cref{def:projprod}. 
\begin{enumerate}
\item {\bf Associativity:} We expand the product as follows:
	\begin{align}
	\begin{split}
	(\bfa \starQ \bfb) \starQ \bfc 
		&= [([(\widehat{\bfa} \odot \widehat{\bfb}) \times_3 \bfQ] \times_3 \bfQ^H) \odot \widehat{\bfc}] \times_3 \bfQ
	\end{split}
	\end{align} 
We then combine the mode-$3$ products by
	\begin{align}
	\begin{split}
	[([(\widehat{\bfa} \odot \widehat{\bfb}) \times_3 \bfQ] \times_3 \bfQ^H) \odot \widehat{\bfc}] \times_3 \bfQ
		&=[[(\widehat{\bfa} \odot \widehat{\bfb}) \times_3 \bfQ^H\bfQ]  \odot \widehat{\bfc}] \times_3 \bfQ\\
		&=[(\widehat{\bfa} \odot \widehat{\bfb})  \odot \widehat{\bfc}] \times_3 \bfQ.
	\end{split}
	\end{align}
Using the associativity of the Hadamard product and reversing the steps completes the proof. 

\item {\bf Distributivity:} We expand the left-hand side as follows:
	\begin{align}
	(\bfa + \bfb) \starQ \bfc 
	&=[ ((\bfa + \bfb) \times_3 \bfQ^H) \odot (\bfc  \times_3 \bfQ^H)] \times_3 \bfQ.
	\end{align}
Using the distributivity of the mode-$3$ and Hadamard products, we get
	\begin{align}
	\begin{split}
	[ ((\bfa + \bfb) \times_3 \bfQ^H) \odot (\bfc  \times_3 \bfQ^H)] \times_3 \bfQ.
	&=[ (\widehat{\bfa} + \widehat{\bfb}) \odot \widehat{\bfc}] \times_3 \bfQ\\
	&=[ \widehat{\bfa} \odot \widehat{\bfc} + \widehat{\bfb} \odot \widehat{\bfc}] \times_3 \bfQ.
	\end{split}
	\end{align}

Using the distributivity of the mode-$3$ product again completes the proof. 
Distribution holds from the left as well following~\Cref{thm:symmetry}. 
	
\end{enumerate}

Because the $\starQ$-product is built on tubal multiplication (see~\eqref{eq:mprodEntrywise}), associativity and distributivity extend to multiplying compatibly-sized tensors as well. 

\end{proof}

\begin{mytheo}{Transposition of the Projected Product}{transpose}
Given $\bfQ\in \Stiefel_{n_3,p}(\Cbb)$ and $\TA\in \Cbb^{n_1\times m\times n_3}$ and $\TB\in \Rbb^{m\times n_2\times n_3}$, we have $(\TA \starQ \TB)^H = \TB^H \starQ \TA^H$. 
\end{mytheo}

\begin{proof}
Let $\bfM = \begin{bmatrix} \bfQ & \bfQ_{\perp}\end{bmatrix}^H$ be a unitary matrix. 
Then, by~\eqref{eq:colSpaceDef}, we have
	\begin{align}
	\TA \starQ \TB =  (\TA \times_3 \bfQ\bfQ^H) \starM  (\TB \times_3 \bfQ\bfQ^H). 
	\end{align}
Using the $\starM$-product definition (\Cref{def:conjugate_transpose}), we have 
	\begin{subequations}
	\begin{align}
	[(\TA \times_3 \bfQ\bfQ^H) \starM (\TB \times_3 \bfQ\bfQ^H)]^H 
	&= (\TB \times_3 \bfQ\bfQ^H)^H \starM  (\TA \times_3 \bfQ\bfQ^H)^H\\
	&= (\TB^H \times_3 \bfQ\bfQ^H) \starM (\TA^H \times_3 \bfQ\bfQ^H) \label{eq:starm_transpose}\\
	&= \TB^H \starQ \TA^H.
	\end{align}
	\end{subequations}
A subtlety in~\eqref{eq:starm_transpose} is that $\TA^H$ and $\TB^H$ are the well-defined $\starM$-conjugate transpose~\cite[Definition 2.1 ]{kilmer2019tensortensor}. 
\end{proof}

A similar theorem can be stated for the $\starQ$-inverse of the $\starQ$-product of two tensors. 
We have demonstrated that the $\starQ$-framework preserves algebraic identities despite the lack of injectivity of certain $\starQ$-operations. 
As a result, we consider the $\starQ$-product is to be matrix mimetic.

%% file: 03_04_generalized_projected_products.tex
\subsection{Generalizing the $\starQ$-product}
\label{sec:starQgeneral}

All definitions and theoretical results extend to nonzero multiples of matrices with unitary columns; that is, for $\bfW = c \bfQ$ where $\bfQ \in \Stiefel_{n_3,p}(\Cbb)$ and $c\in \Cbb\backslash \{0\}$ .  
In this case, we would define $\starW$-product using the pseudoinverse and obtain the following relationship:
	\begin{subequations}
	\begin{align}
	\TA \starW \TB
	 &=[(\TA\times_3 \bfW^H) \smalltriangleup (\TB \times_3 \bfW^H)] \times_3 (\bfW^H)^{\dagger}\\
	&=[(\TA\times_3 (c\bfW^H) \smalltriangleup (\TB \times_3 (c\bfW)^H)] \times_3 (\tfrac{1}{c}\bfW)\\
	&= c (\TA \starQ \TB).
	\end{align} 
	\end{subequations}
Because $\TA \star_{\bfW^H}' \TB$ is a scalar multiple of $\TA \starQ \TB$, linear algebraic properties will be preserved. For ease of presentation and discussion, we provide properties and theory for the $\starQ$-product with the understanding that these properties extend to the $\star_{\bfW^H}'$-product as well.

%% file: 04_eckart_young.tex
\section{The $\starQ$-SVD and Eckart-Young Theorem}
\label{sec:tsvdq}

With the algebraic building blocks in place, we have the tools to define a $\starQ$-based tensor SVD, which strongly resembles the $t$-SVDM in~\cite{kilmer2019tensortensor}. We prove the optimality of low-rank representations in~\Cref{thm:projEckartYoung} and provide insight into an optimal choice of transformation matrix in~\Cref{thm:projection_error}. In~\Cref{sec:tsvdqII}, we present a more compressible variant of the $\starQ$-SVD and related theoretical results. In~\Cref{sec:hosvd}, we compare the $\starQ$-SVD to the truncated higher-order SVD~\cite{LathauwerMoorVandewalle2000}. 

\begin{definition}[$\starQ$-SVD]\label{def:projtsvdq}

Given $\bfQ\in \Stiefel_{n_3,p}(\Cbb)$ and tensor $\TA\in \Cbb^{n_1\times n_2\times n_3}$, the $\starQ$-SVD is 
	\begin{align}
	\TA \times \bfQ \bfQ^H =  \TU \starQ \TS \starQ \TV^H
	\end{align}
where $\TU\in \Cbb^{n_1\times n_1\times n_3}$ and $\TV\in \Cbb^{n_2\times n_2\times n_3}$ are $\starQ$-unitary and $\TS\in \Cbb^{n_1\times n_2\times n_3}$ is f-diagonal with
	\begin{align}
	\|\TS_{1,1,:} \|_F \ge \|\TS_{2,2,:}\|_F \cdots \ge \|\TS_{q,q,:}\|_F \ge 0 \qquad \text{for $q = \min(n_1,n_2)$.}
	\end{align}
\end{definition}

We present the pseudocode to compute the (truncated) $\starQ$-SVD, including the computational costs for dense matrix operations from~\cite[Section 1.4.1]{Halko2011}. Internally, the algorithm relies on the matrix SVD, hence the $\starQ$-SVD exists for all tensors. 
\begin{myalg}{Truncated $\starQ$-SVD}{tsvdq}

\begin{algorithmic}[1]
\State \textbf{Inputs:} $\TA\in \Cbb^{n_1\times n_2\times n_3}$, $\bfQ\in \Stiefel_{n_3,p}(\Cbb)$, truncation parameter $k \in \{1,\dots,\min(n_1,n_2)\}$
 
\State Move to the transform domain $\widehat{\TA} = \TA \times_3 \bfQ^H$ \Comment{$\Ocal(n_1n_2n_3p)$}
\State Compute truncated matrix SVDs of frontal slices \Comment{$\Ocal(n_1n_2k p)$}
    
    	$$\widehat{\TA}_{:,:,i}  \approx \widehat{\TU}_{:,1:k,i} \widehat{\TS}_{1:k,1:k,i} \widehat{\TV}_{:,1:k,i}^H \qquad \text{for $i=1,\dots,p$.}$$

    \State Return to the spatial domain  \Comment{$\Ocal((n_1k + k^2 + kn_2)n_3p)$}
    	\begin{align*}
	\TU_k = \widehat{\TU}_k \times_3 \bfQ, \quad
	\TS_k = \widehat{\TS}_k \times_3 \bfQ,\quad \text{and}\quad
	\TV_k = \widehat{\TV}_k \times_3 \bfQ.
	\end{align*}
	
\State \textbf{Return:} $\TU_k, \TS_k, \TV_k$

 \end{algorithmic}
 
\end{myalg}

The $\starQ$-SVD is less expensive computationally and storage-wise than the original $\starM$-SVD (or $t$-SVDM in~\cite[Algorithm 2]{kilmer2019tensortensor}) by a factor of $p / n_3$; see~\Cref{fig:starM_vs_starQ} for intuition.  

The $\starQ$-SVD gives rise to a notion of tensor rank consistent with that of matrix rank, traditionally called $t$-rank if the $\starM$-literature. 

\begin{definition}[$\starQ$-rank]\label{def:starMRank}
Given a tensor and its $\starQ$-SVD $\TA \approx \TU \starQ \TS \starQ \TV^H$, the $\starQ$-rank of $\TA$ is the number of nonzero tubes in $\TS$; that is, 
   	\begin{align}
	\starQ\text{-rank}(\TA) = \#\{i \mid \|\TS_{i,i,:}\|_F > 0\text{ for }i=1,\dots,n_3\}
	\end{align}
where $\#$ denotes the cardinality of the set. 
\end{definition}

Using $\bfM$ from~\Cref{sec:equivalent_presentations} the $\starM$-SVD from~\cite{kilmer2019tensortensor} and the $\starQ$-SVD are share information. Specifically, in the transform domain, the $\starQ$-SVD factors are equal to the first $p$ frontal slice factors of the $\starM$-SVD. We can thus connect the two notions of $\starQ$-rank and $\starM$-rank.

\begin{mytheo}{$\starQ$-rank  $\le$ $\starM$-rank}{trank}
Let $\bfM\in \Stiefel_{n_3,n_3}(\Cbb)$ and $\bfQ = \bfM_{1:p,:}^H$. Then, $\projtrank(\TA) \le \trank(\TA)$ for any tensor $\TA$. 
\end{mytheo}

\begin{proof}
An equivalent definition of $\starQ$-rank (\Cref{def:starMRank}) is the maximum rank of the frontal slices in the transform domain; that is, 
	\begin{subequations}
	\begin{align}
	\trank(\TA) 	&= \max_{i=1,\dots,n_3} \rank(\TA \times_3 \bfM_{i,:})\\
	\projtrank(\TA) 	&= \max_{i=1,\dots,p} \rank(\TA \times_3 \bfQ_{:,i}^H)
	\end{align}
	\end{subequations}
Because $\bfQ_{:,i}^H = \bfM_{i,:}$ for $i=1,\dots, p$, the first $p$ frontal slices of $\TA \times_3 \bfM$ are equal to the $p$ frontal slices of $\TA \times_3 \bfQ^H$. Thus, $\projtrank(\TA)$ can be no larger than $\trank(\TA)$. 

\end{proof}

While the $\starQ$-SVD is not equal to the original tensor; i.e., $\TA \not= \TU \starQ \TS \starQ \TV^H$, it still 
satisfies an Eckart-Young-like optimality result. 

\begin{mytheo}{$\starQ$-Eckart-Young Optimality}{projEckartYoung}

Given $\bfQ \in \Stiefel_{n_3,p}(\Cbb)$ and $\starQ$-SVD $\TA \approx \TU \starQ \TS \starQ \TV^H$ with $\projtrank(\TA) = r$, the best $\starQ$-rank-$k$ approximation to $\TA$ for $k\le r$ is given by the truncated $\starQ$-SVD; that is,
	\begin{align}
	\TA_k \equiv \TU_{:,1:k,:} \starQ \TS_{1:k,1:k,:} \starQ \TV_{:,1:k,:}^H \in \argmin_{\TB \in \Bcal_k'} \|\TA - \TB\|_F
	\end{align}
where $\Bcal_k' = \{\TX \starQ \TY^H \mid \TX\in \Cbb^{n_1\times k\times n_3}, \TY\in \Cbb^{n_2\times k\times n_3}\}$.  The Frobenius norm error is
	\begin{align}
	\|\TA - \TA_k\|_F^2 = \underbrace{\sum_{j=k+1}^{r} \|\TS_{j,j,:}\|_F^2}_{\text{Eckart-Young error}} + \underbrace{\|\TA \times_3 (\bfI_{n_3} - \bfQ \bfQ^H)\|_F^2 \vphantom{\sum_{j=k+1}^{r} \|\TS_{j,j,:}\|_F^2}}_{\text{projection error}}.
	\end{align}
\end{mytheo}

\begin{proof} 
Recall from~\eqref{eq:tensorSplitProjection}, we can separate the tensor as $\TA = \TA \times_3 \bfQ \bfQ^H + \TA \times_3 (\bfI_{n_3} - \bfQ \bfQ^H)$. Following~\Cref{def:projtsvdq}, the $\starQ$-SVD can exactly capture the first term; that is, 
	\begin{align}
	\TA \times_3 \bfQ \bfQ^H  = \TU \starQ \TS \starQ \TV^H.
	\end{align}
The second term lies in the null space of $\bfQ^H$, and hence cannot be approximated by the full nor truncated $\starQ$-SVD. Thus, we quantify the optimal $\starQ$-rank-$k$ approximation of $\TA \times_3 \bfQ \bfQ^H$. 

Using the unitary invariance of the Frobenius norm, we have 
	\begin{align}
	\|\TA\times_3 \bfQ \bfQ^H - \TX \starQ \TY^H\|_F = \|\widehat{\TA}- \widehat{\TX} \smalltriangleup \widehat{\TY}^H\|_F 
	\end{align}
where $\widehat{\cdot} = \cdot \times_3 \bfQ^H$. In the transform domain, we independently approximate each frontal slice. By the matrix Eckart-Young Theorem~\cite{EckartYoung}, we have
	\begin{align}
	\|\widehat{\TA}_{:,:,i} - \widehat{\TA}_k(:,:,i) \|_F \le \|\widehat{\TA}_{:,:,i} - \widehat{\TX}_{:,:,i} \widehat{\TY}_{:,:,i}^H\|_F
	\end{align}
for $i=1,\dots,p$ where $\TA_k$ is the truncated $\starQ$-SVD. This completes the proof that the truncated $\starQ$-SVD produces the best $\starQ$-rank-$k$ approximation to the original tensor.  

To compute the two error terms, let $\bfM = \begin{bmatrix} \bfQ & \bfQ_{\perp}\end{bmatrix}^H$ be a unitary matrix. By~\eqref{eq:backFrontalSlicesZero}, the last $n_3 - p$ frontal slices of $\TA_k \times_3 \bfM$ are zero.  By the unitary invariance of the Frobenius norm, 
	\begin{align}
	\|\TA - \TA_k\|_F^2 
		=\|(\TA - \TA_k) \times_3 \bfM\|_F^2 
		=\|\TA \times_3 \bfQ^H - \TA_k \times_3 \bfQ^H\|_F^2 + \|\TA \times_3 \bfQ_{\perp}^H - {\bf0}\|_F^2
	\end{align}
The first term is the error of the $\starQ$-SVD approximation of the first $p$ frontal slices in the transform domain (Eckart-Young error). Following~\cite[Theorem 3.7]{kilmer2019tensortensor}, the Eckart-Young error is equal to the norm of the truncated singular tubes. The second term is equal to $\|\TA \times_3  (\bfI_{n_3} - \bfQ \bfQ^H)\|_F^2$ because $\bfQ_{\perp}\bfQ_{\perp}^H = \bfI_{n_3} - \bfQ \bfQ^H$ (projection error). 

\end{proof}

\begin{myrem}
A more descriptive notation for the truncated $\starQ$-SVD would be $\TA_{k,p}$. However, because the $\starQ$-dependence implies the projection dimension $p$, we prefer the more conventional $\TA_k$ notation and assume an implicit dependence on $\bfQ$ and $p$.
\end{myrem}

For a specific choice of transformation matrix, we can specify a concrete projection error. 

\begin{figure}
    \centering
    \includegraphics[width=1\linewidth]{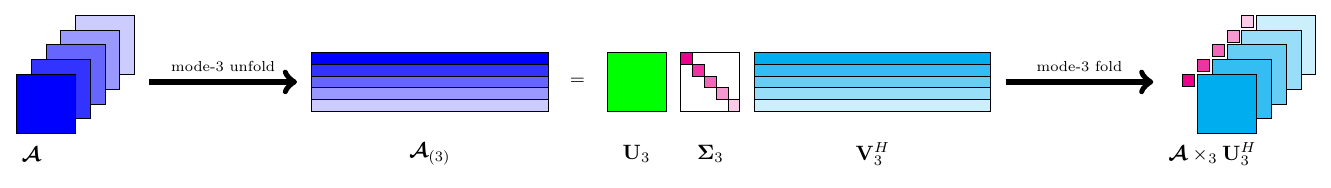}
    \caption{Illustration of $\TA \times_3 \bfU_3^H$ where $\bfU_3$ is the left-singular matrix from the mode-$3$ unfolding of $\TA$; that is, $\bfA_{(3)} = \bfU_3 \bfSigma_3 \bfV_3^H$. The frontal slices of $\TA \times_3 \bfU_3^H$ have Frobenius norm equal to the singular values of the mode-$3$ unfolding. From the ordering of the singular values (indicated by the various shades of magenta), the relative importance of each transformed frontal slice decays from front (dark magenta) to back (light magenta).}
    \label{fig:mode3svd}
\end{figure}

\begin{mycorollary}{$\starQ$-SVD Error for $\bfU_3$}{projEckartYoungZ}
Let $\bfU_3$ be the left-singular vectors of the mode-$3$ unfolding of a tensor $\TA$; that is, $\bfA_{(3)} = \bfU_3 \bfSigma_3 \bfV_3^H$. Then, with $\bfQ = \bfU_3(:,1:p)$, the truncated $\starQ$-SVD error is 
	\begin{align}
	\| \TA - \TA_k\|_F^2 = \sum_{j=k+1}^r \|\TS_{j,j,:}\|_F^2 + \sum_{j=p+1}^{n_3} \sigma_j(\bfA_{(3)})^2
	\end{align}
where $\sigma_j(\bfY)$ returns the $j$-th largest singular value of a matrix $\bfY$. 
\end{mycorollary}

\begin{proof} 
The proof follows directly from~\Cref{thm:projEckartYoung}.  The second term comes from the structure of the transformation matrix $\bfU_3^H$ (see~\Cref{fig:mode3svd}). Specifically,  in the transform domain, we have
	\begin{align}
	(\TA \times_3 \bfU_3^H)_{:,:,i} = \sigma_i(\bfA_{(3)}) \texttt{reshape}(\bfV_3(:,i),[n_1,n_2])
	\end{align}
for $i=1,\dots,n_3$. The $\starQ$-SVD cannot approximation the $p+1$ through $n_3$ frontal slices in the transform domain. Hence, the error of those frontal slices is
 \begin{align}
 \|{\bf0} -\sigma_i(\TA_{(3)}) \texttt{reshape}(\bfV_3(:,i),[n_1,n_2]) \|_F = \sigma_i(\TA_{(3)}). 
 \end{align}
Because the columns of $\bfV_3$ have unit length, the magnitude depends only on the singular value. 
\end{proof}

The $\starQ$-SVD with $\bfQ = \bfU_3(:,1:p)$ optimizes the projection error in the following sense.

\begin{mytheo}{Optimal Projection Error}{projection_error}
Given $\bfA_{(3)} = \bfU_3\bfSigma_3 \bfV_3^H$, then $\bfU_3(:,1:p) \in \argmin_{\bfQ\in \Stiefel_{n_3,p}(\Cbb)} \|\TA \times_3 (\bfI_{n_3} - \bfQ\bfQ^H)\|_F.$
\end{mytheo}

\begin{proof} 
We can write the projection error as
\begin{align}
\|\TA \times_3 (\bfI_{n_3} - \bfQ\bfQ^H)\|_F 
	&=\| (\bfI_{n_3} - \bfQ\bfQ^H) \bfA_{(3)}\|_F =\| \bfA_{(3)} - \bfQ\bfQ^H \bfA_{(3)}\|_F
\end{align}
Note that  $\rank(\bfQ\bfQ^H \bfA_{(3)}) \le p$.  By the matrix Eckart-Young Theorem~\cite{EckartYoung}, the best rank-$p$ approximation to $\bfA_{(3)}$ in the Frobenius norm is the truncated matrix SVD.  This corresponds to the case when $\bfQ = \bfU_3(:,1:p)$.  
\end{proof}

While the matrix $\bfQ = \bfU_3(:,1:p)$ produces the smallest projection error, it does not necessarily yield the smallest $\starQ$-SVD error (\Cref{thm:projEckartYoung}).  
We provide a counterexample in~\Cref{exam:tsvdmq_error_counter} to illustrate this point. 

\begin{myexam}{Counterexample of Optimal $\starQ$-SVD Error (\Cref{thm:projection_error})}{tsvdmq_error_counter}

Consider the following $2\times 2\times 2$ tensor
	\begin{align}
	\TA_{:,:,1} = \begin{bmatrix} 1 & 0 \\ 0 & 1 \end{bmatrix}  
	\qquad \text{and} \qquad
	\TA_{:,:,2} = \begin{bmatrix} 1 & 0 \\ 0 & -1\end{bmatrix}.
	\end{align}
Then, up to column permutation and negation,  $\bfU_3 = \bfI_2$. 
We compare the truncated $\starQ$-SVD for $\bfU_3$ and the transposed $2\times 2$ Haar wavelet matrix 
	\begin{align}
	\bfH_2 = \tfrac{1}{\sqrt{2}}\begin{bmatrix} 1 & 1 \\ 1 & -1 \end{bmatrix}.
	\end{align}
Using truncation parameter $k=1$ and projection dimension $p=1$, we compute the approximation errors
	\begin{align}
	 \|\TA - \TA_1(\bfU_3(:,1))\|_F^2 &= 1 + 2 \qquad \text{and}\\
	  \|\TA - \TA_1(\bfH_2(1,:)^\top)\|_F^2 &= 0 + 2
	\end{align}
where $\TA_1(\bfQ)$ is the projected $\starQ$-SVD obtained using matrix $\bfQ$. 
Here, $2$ is the squared projection error obtained by not approximating the second frontal slice in the transform domain. 
For $\bfQ = \bfH(1,:)^\top$, the transform domain frontal slice is
	\begin{align}
	\TA \times_3 \bfH_2(1,:) = \begin{bmatrix} \sqrt{2} & 0 \\ 0 & 0 \end{bmatrix} 
	\end{align} 
Since this matrix is rank-$1$, it can be approximated exactly by the rank-$1$ truncated matrix SVD. 
This means the Eckart-Young error is zero for this matrix. 
\end{myexam}

%% file: 04_01_tsvdqII.tex
\subsection{$\starQ$-SVDII}
\label{sec:tsvdqII}

In~\cite{kilmer2019tensortensor}, a variant called the $\starM$-SVDII (or often the $t$-SVDMII), was proven to have the same approximation quality as the $\starM$-SVD for less storage cost. We present a similar variant for projected products called the $\starQ$-SVDII.  The key to the additional compressibility is to consider a global perspective of the $\starQ$-SVD. 
By viewing the facewise product as block diagonal matrix multiplication, the $\starQ$-SVD in the transform domain can be written as 
	\begin{align}\label{eq:blkdiag}
	\setlength\arraycolsep{1pt}
	\begin{bmatrix}
	\widehat{\TA}_{:,:,1}\\
	&\ddots \\
	&&\widehat{\TA}_{:,:,p}
	\end{bmatrix} &=
	\setlength\arraycolsep{1pt}
	 \begin{bmatrix}
	\widehat{\TU}_{:,:,1}\\
	&\ddots \\
	&&\widehat{\TU}_{:,:,p}
	\end{bmatrix}
	\setlength\arraycolsep{1pt}
	 \begin{bmatrix}
	\widehat{\TS}_{:,:,1}\\
	&\ddots \\
	&&\widehat{\TS}_{:,:,p}
	\end{bmatrix}
	\setlength\arraycolsep{1pt}
	 \begin{bmatrix}
	\widehat{\TV}_{:,:,1}\\
	&\ddots \\
	&&\widehat{\TV}_{:,:,p}
	\end{bmatrix}^H.
	\end{align}

From the block diagonal presentation in~\eqref{eq:blkdiag}, we can globally reorder the singular values and  truncate based on a desired approximation quality.   
Let $r$ be the $\starQ$-rank of $\TA$ and let $\pi:  \{(j,j,i) | j=1,\dots,r,  i=1,\dots,p\} \to \{1,\dots, rp\}$ be a permutation that sorts the transform domain singular values in decreasing order of magnitude; i.e., 
	\begin{align}\label{eq:singvalreorder}
	\widehat{\bfs}_{\pi(j,j,i)} = \widehat{\TS}_{j,j,i} \quad \text{for }j=1,\dots,r \text{ and } i=1,\dots,p
	\end{align}
such that $\widehat{\bfs}_{\ell} \ge \widehat{\bfs}_{\ell+1}$ for $\ell=1,\dots,rp-1$.  
Note that the mapping $\pi$ is not unique if there are repeated singular values. 
Using this reordering, we truncate the $\starQ$-SVDII based on a user-defined energy level $\gamma \in (0,1]$.  
We present the pseudocode in~\Cref{alg:tsvdqII_multirank}

\begin{myalg}{$\starQ$-SVDII}{tsvdqII_multirank}

\begin{algorithmic}[1]
\State \textbf{Inputs:} $\TA\in \Cbb^{n_1\times n_2\times n_3}$, $\bfQ\in \Stiefel_{n_3,p}(\Cbb)$, energy parameter $\gamma \in (0,1]$

\State Move to the transform domain $\widehat{\TA} = \TA \times_3 \bfQ^H$
\State Compute the facewise SVD $\TA = \widehat{\TU} \smalltriangleup \widehat{\TS} \smalltriangleup \widehat{\TV}^H$

\State Globally reorder singular values and store in a vector $\widehat{\bfs}$ according to~\eqref{eq:singvalreorder}  \label{alg:energy_truncation}

\State Find the first index $K$ such that {$\sfrac{\sum_{\ell=1}^{K}\widehat{\bfs}_\ell^2}{\|\widehat{\bfs}\|_2^2} \ge \gamma$}

\For{$i=1,\dots,p$}
\State Set the rank per frontal slice to be $\rho_i = \argmax_{j=1,\dots,r} \{\pi(j,j,i) \le K\}$
\State Truncate in the transform domain with
	\begin{align*}
	\widehat{\TU}_{:,:,i} &\gets \widehat{\TU}_{:,1:\rho_i,i}, &
	\widehat{\TS}_{:,:,i} &\gets \widehat{\TS}_{1:\rho_i,1:\rho_i,i}, & \text{and}&&
	\widehat{\TV}_{:,:,i} &\gets \widehat{\TV}_{:,1:\rho_i,i}.
	\end{align*}
\EndFor

\State {\bf Return:} $\widehat{\TU}, \widehat{\TS}, \widehat{\TV}$

 \end{algorithmic}
 
 \end{myalg}

A key difference between the $\starQ$-SVDII and the $\starM$-SVDII~\cite[Algorithm 3]{kilmer2019tensortensor} is the total energy in the transform domain may be different than total energy in the spatial domain; that is, $\|\TA \times_3 \bfQ^H\|_F \not= \|\TA\|_F$.   
In particular, the denominator in Line 5 of~\Cref{alg:tsvdqII_multirank} has an implicit dependence on $\bfQ$. 

The $\starQ$-SVDII gives rise to two different notions of the rank of a tensor.

\begin{definition}[$\starQ$-multirank]\label{def:multirank}
The $\starQ$-multirank of $\TA$ is a $p$-tuple $\bfrho$ where $\rho_i = \rank(\widehat{\TA}_{:,:,i})$ for $i=1,\dots,p$ and $\widehat{\TA} = \TA \times_3 \bfQ^H$. 
\end{definition}

\begin{definition}[$\starQ$-implicit rank]\label{def:implicit_rank}
The $\starQ$-implicit rank of $\TA$ is the total number of singular values stored in the transform domain; that is,  $\rho=\sum_{i=1}^p \bfrho_i$ where $\bfrho$ is the $\starQ$-multirank of $\TA$. 
\end{definition}
If a tensor has $\starQ$-rank-$r$, then the tensor has a $\starQ$-implicit rank of $\rho \le rp$. 
Combining these notions of rank, we prove the $\starQ$-SVDII is provably optimal in an Eckart-Young sense.  

\begin{mytheo}{Optimality of $\starQ$-SVDII}{tsvdqII} 

Given $\bfQ \in \Stiefel_{n_3,p}(\Cbb)$ and a tensor $\TA\in \Cbb^{n_1\times n_2\times n_3}$ with $\starQ$-implicit rank $\rho$,   the best $\starQ$-multirank-$\bfkappa$ approximation $\TA_{\bfkappa}$ with $\sum_{i=1}^p \bfkappa_i = \kappa \le \rho$ is the $\starQ$-SVDII; that is,
	\begin{align}
	\TA_{\bfkappa} \in \argmin_{\TB \in \Bcal_k''} \|\TA - \TB\|_F
	\end{align}
where $\Bcal_k'' = \{\TX \in  \Cbb^{n_1\times n_1 \times n_3} \mid \starQ\text{-implicit-rank}(\TX) \le \kappa\}$.  
The Frobenius norm error is
	\begin{align}
	\|\TA - \TA_{\bfkappa}\|_F^2 = \sum_{i=1}^{p} \sum_{j=\bfkappa_i + 1}^{r} \widehat{\TS}_{j,j,i}^2 + \|\TA \times_3 (\bfI_{n_3} - \bfQ \bfQ^H)\|_F^2
	\end{align}
\end{mytheo}

\begin{proof}
The proof leverages the matrix Eckart-Young Theorem similarly to the proof of~\Cref{thm:projEckartYoung}, hence we omit the details for brevity. 
\end{proof}

We further prove that the $\starQ$-SVDII computed via~\Cref{alg:tsvdqII_multirank} offers a better approximation than the $\starQ$-SVD for no additional storage cost. 

\begin{mytheo}{$\starQ$-SVDII vs. $\starQ$-SVD}{tsvdqII_vs_tsvdq}
Given $\bfQ\in \Stiefel_{n_3,p}(\Cbb)$ and $\TA \in \Cbb^{n_1\times n_2\times n_3}$ with $\starQ$-rank-$r$, let $\TA_k$ of the truncated $\starQ$-SVD for $k \le r$.  
Define the energy parameter $\gamma = \sfrac{\|\TA_k \times_3 \bfQ^H\|_F^2}{\|\TA \times_3 \bfQ^H\|_F^2}$. 
Then, the $\starQ$-SVDII $\TA_{\bfkappa}$ corresponding to energy $\gamma$ is a no worse approximation than $\TA_k$; that is, 
	\begin{align}
	\|\TA - \TA_{\bfkappa}\|_F \le \|\TA - \TA_k\|_F.
	\end{align}
\end{mytheo}

\begin{proof} 
Because we use the same transformation $\bfQ$ for both decompositions, the projection error is equal for the $\starQ$-SVD and $\starQ$-SVDII. 
Thus, we will compare that the Eckart-Young error terms in~\Cref{thm:projEckartYoung} and~\Cref{thm:tsvdqII}. 
For the $\starQ$-SVDII, if we use the trivial $\starQ$-multirank is $\bfkappa = (k,k,\dots,k)$, then $\TA_{\bfkappa} = \TA_k$ and the relative energy constraint is automatically satisfied.  
The $\starQ$-multirank formed from the energy threshold in~\Cref{alg:tsvdqII_multirank} will truncate the smallest singular values globally.  Thus, we can only improve the approximation compared to the trivial $\starQ$-multirank. 
As a result, the Eckart-Young (EY) error term satisfies the following inequality:
	\begin{align}
	\underbrace{\sum_{i=1}^p \sum_{i=\rho_i+1}^r\widehat{\TS}_{j,j,:}^2}_{\text{\begin{tabular}{c} $\starQ$-SVDII EY error\\ with energy $\gamma$ \end{tabular}}} \le 
	\underbrace{\sum_{i=1}^p \sum_{i=k+1}^r\widehat{\TS}_{j,j,:}^2}_{\text{\begin{tabular}{c}$\starQ$-SVDII EY error\\ with $\bfrho=(r,\dots,r)$ \end{tabular}}}
	 =  \sum_{j=k+1}^r \|\widehat{\TS}_{j,j,:}\|_F^2
	= \underbrace{\sum_{j=k+1}^r \|\TS_{j,j,:}\|_F^2.}_{\text{$\starQ$-SVD EY error}}
	\end{align}

\end{proof}

%% file: 04_02_hosvd.tex
\subsection{Comparison to Higher-Order SVD}
\label{sec:hosvd}

In addition to comparing the $\starQ$-SVD to the original $\starM$-version, we compare to the commonly-used (truncated) higher-order SVD (HOSVD)~\cite{LathauwerMoorVandewalle2000}, which approximates a third-order tensor $\TA\in \Cbb^{n_1\times n_2\times n_3}$ as
	\begin{align}
	\TA \approx \TA_{\bfk} = \TG \times_1 \bfU_1(:,1:k_1) \times_2 \bfU_2(:,1:k_2) \times_3 \bfU_3(:,1:k_3)
	\end{align}
where $\TG\in \Rbb^{k_1\times k_2\times k_3}$ is the core tensor and $\bfU_i \in \Stiefel_{n_i,n_i}(\Cbb)$ for $i=1,\dots,3$ are the factor matrices. 
Each factor matrix $\bfU_i$ is the truncated left singular matrix from the SVD of various tensor unfoldings
	\begin{align}\label{eq:modek_svd}
	\bfA_{(i)} = \bfU_i \bfSigma_i \bfV_i^H \qquad \text{for $i=1,2,3$}
	\end{align}
where $\bfA_{(i)}$ is the mode-$i$ unfolding, defined similarly to~\Cref{def:mode3unfolding}; details can be found in~\cite{KoldaBader}.  We denote the truncated HOSVD with multilinear rank $\bfk = (k_1,k_2,k_3)$ as $\TA_{\bfk}$. We prove that the truncated $\starQ$-SVD yields a more accurate approximation than the truncated HOSVD. 

\begin{mytheo}{$\starQ$-SVD vs. HOSVD}{projhosvd} 
If $\bfQ =\bfU_3(:,1:p)$ and $k=\min(k_1,k_2)$, then the $\starQ$-rank-$k$ approximation $\TA_k$ is no worse than the truncated HOSVD $\TA_{\bfk}$ of multilinear rank $\bfk = (k_1,k_2,p)$; i.e., 
	\begin{align}
	\|\TA - \TA_k \|_F \le \|\TA - \TA_{\bfk}\|_F.
	\end{align}
\end{mytheo}

\begin{proof} 
The key observation is that the transformed tensor, $\widehat{\TA} = \TA \times_3 \bfU_3(:,1:p)^H$, is found in the HOSVD as well. Specifically, we regroup the truncated HOSVD as
	\begin{subequations}
	\begin{align}
	\TG  &= (\TA \times_3 \bfU_3(:,1:p)^H) \times_1 \bfU_1(:,1:k_1)^H \times_2 \bfU_2(:,1:k_2)^H\\
		&=\widehat{\TA} \times_1 \bfU_1(:,1:k_1)^H \times_2 \bfU_2(:,1:k_2)^H
	\end{align}
	\end{subequations}
The remainder of the proof follows from~\cite[Section 6.A and Theorem 6.1]{kilmer2019tensortensor}, which shows that each frontal slice of the truncated HOSVD has rank less than or equal to $k$. 
Because the $\starQ$-SVD uses the optimal rank-$k$ approximation to each frontal slice, it yields a more accurate approximation. 
\end{proof}

Combining~\Cref{thm:projhosvd} and~\Cref{thm:tsvdqII_vs_tsvdq}, we can produce a similar theorem for the $\starQ$-SVDII. 
We omit the theorem and proof for the sake of brevity.

%% file: 05_numerical_experiments.tex
\section{Numerical Experiments}
\label{sec:numericalexperiments}
We present several numerical experiments to demonstrate the approximation quality and compressibility of the $\starQ$-SVD and variants. We compare $\starQ$-SVD for various choices of transformation $\bfQ$ on two gray-scale videos datasets (\Cref{sec:video}) and on a hyperspectral imaging dataset (\Cref{sec:hyperspectral}). In~\Cref{sec:tsvdq_vs_hosvd}, we compare the $\starQ$-SVDII to the truncated HOSVD  on  the same hyperspectral imaging dataset, providing empirical support of the theory presented in~\Cref{sec:hosvd}. All code and experiments are available at \url{https://github.com/elizabethnewman/projected-products.git}.

\subsection{Experiment Setup}

Throughout the presented experiments, we will common transformation matrices and approximation and compression metrics. 
We present the details here for concision. 

\subsubsection{Transformation Matrices}

We compare four transformation matrices in our experiments: the identity matrix $\bfI$, a random orthogonal matrix\footnote{In {\sc Matlab}, we write \texttt{W = orth(randn(n3)); Q = W(:,1:p);}} $\bfW$, the (tranposed) discrete cosine transform (DCT) matrix\footnote{In {\sc Matlab}, we write \texttt{C = dctmtx(n3); Q = C(1:p,:)';}}  $\bfC^\top$, and the data-dependent, left-singular vectors of the mode-$3$ unfolding $\bfU_3$. The identity and random matrices are control cases; the identity does not exploit any correlations along the third-dimension, but does not require any additional storage; the random matrix has no prescribed structure. The DCT and data-dependent matrices are structured cases; the DCT matrix can be interpreted as a real-valued approximation of the $t$-product and the data-dependent matrix produces the optimal projection error (\Cref{thm:projection_error}). We only consider real-valued transformations to avoid introducing complex values in the approximations.

\subsubsection{Metrics}

We consider two metrics: relative error (RE) to measure approximation quality and compression ratio (CR) to measure storage costs. We define each metric as
	\begin{align}
	\text{RE} = \frac{\|\TA - \widetilde{\TA}\|_F}{\|\TA\|_F} \qquad \text{and}
	\qquad \text{CR} = \frac{\texttt{st}[\TA]}{\texttt{st}[\widetilde{\TA}] + \texttt{st}[\bfQ]}
	\end{align}
where $\widetilde{\TA}$ is the compressed representation of the original data $\TA$ and $\texttt{st}[\cdot]$ computes the storage cost of the input. We seek small relative errors (close to zero) and large compression ratios (greater than one indicates that we have achieved compression). 

In our experiments, we assume we are given a dense, real-valued tensor $\TA\in \Rbb^{n_1\times n_2\times n_3}$ with storage cost $\texttt{st}[\TA] = n_1 n_2 n_3$.  For a truncation parameter $k \le \min(n_1,n_2)$, the storage cost for the $\starQ$-SVD approximation $\TA_k$ is 
	\begin{align}
	\text{$\starQ$-SVD:} & \quad \texttt{st}[\TA_k]  = \texttt{st}[\widehat{\TU}_k] + \texttt{st}[\widehat{\TS_k \starQ \TV_k^H}]  = n_1 k p + k n_2 p.
	\end{align}
For an implicit rank of $\kappa \le \min(n_1,n_2) p$, the storage cost for the $\starQ$-SVDII approximation $\TA_{\bfkappa}$ is  
	\begin{align}
	\text{$\starQ$-SVDII:} & \quad \texttt{st}[\TA_{\bfkappa}]  = \texttt{st}[\widehat{\TU}_{\bfkappa}] + \texttt{st}[\widehat{\TS_{\bfkappa} \starQ \TV_{\bfkappa}^H}]= \kappa (n_1 + n_2).
	\end{align}
Note that we store the $\starQ$-approximations in the transform domain.  In general, we require additional storage of the transformation matrix $\texttt{st}[\bfQ] = n_3 p$, except when we use an identity transformation and have no additional overhead.  

The truncated HOSVD with multirank $\bfk = (k_1,k_2,k_3)$ has a storage cost of
	\begin{align}
		\text{HOSVD:} & \quad \texttt{st}[\TA_{\bfk}]  = \texttt{st}[\TG] + \texttt{st}[\bfU_1] + \texttt{st}[\bfU_2]  + \texttt{st}[\bfU_3] = k_1 k_2 k_3 + n_1 k_1 + n_2 k_2 + n_3 k_3.
	\end{align}

We also compare performance of the $\starQ$-approximations to the equivalent matrix SVD approximations. Following~\cite[Theorem 5.3]{kilmer2019tensortensor}, we reshape the tensor into a matrix by stacking the frontal slices vertically; that is, $\bfA = \texttt{unfold}[\TA] = \begin{bmatrix} \TA_{:,:,1}^H & \cdots & \TA_{:,:,n_3}^H \end{bmatrix}^H$. The resulting matrix $\bfA$ is of size $n_1n_3\times n_2$. We then compute the matrix SVD $\bfA = \bfU \bfSigma \bfV^H$, which has a storage cost of 
	\begin{align}
	\text{matrix SVD:} & \quad \texttt{st}[\bfA] =  \texttt{st}[\bfU] +  \texttt{st}[\bfSigma \bfV^H] =  n_1n_3k + kn_2. 
	\end{align}

%% file: 05_01_video.tex
\subsection{Video Compression}
\label{sec:video}
We illustrate the utility of the $\starQ$-representations for video compression on two datasets, \texttt{traffic} and \texttt{shuttle}, included in the {\sc Matlab} Image Processing Toolbox (see~\Cref{sec:video_data} for a visualization). Both datasets are oriented as $\text{height} \times \text{width} \times \text{time}$; i.e., if $\TA$ is a video, then $\TA_{:,:,i}$ is a grayscale image representing the $i$-th frame. The \texttt{traffic} video ($120 \times 160 \times 120$) consists of a static background (road) and a dynamic foreground (moving cars) recorded from a fixed camera location. The \texttt{shuttle} video ($288 \times 512 \times 121$) captures a rocket launch with a camera following the rocket's vertical trajectory and the rocket exhaust changing the background. We examine the relative error of the low-rank $\starQ$-SVD approximations for the \texttt{traffic} video in~\Cref{fig:traffic} and for the \texttt{shuttle} video in~\Cref{fig:shuttle} for various choices of truncation parameters $k$, projection dimensions $p$, and transformations $\bfQ$.  Because of the similarities in behavior, we create one figure per video and discuss the commonalities and differences in performance subsequently. 

\def\myname{traffic}

\input{05_01_video_figure}

\def\myname{shuttle}
\input{05_01_video_figure}

In~\Cref{fig:traffic_k5} and~\Cref{fig:shuttle_k5} for $k=5$, we observe that the transformations that illuminate practical data structure, discrete cosine and data-dependent transforms, are able to achieve the optimal approximation with only about $10\%$ of the frontal slices stored in the transform domain (i.e., $p/n_3 \approx 0.1$). In comparison, the transformations that do not exploit multilinear structure well, the identity and random matrices, obtain relative errors about an order of magnitude larger for all choices of $p < n_3$ and to not reach the optimal relative error until $p = n_3$.

This pattern of approximation quality is not unique behavior for $k=5$.  Qualitatively, in~\Cref{fig:traffic_approx} and~\Cref{fig:shuttle_approx}, we observe that  the quality for the identity and random matrices improves abruptly when enough information is retained whereas the approximations using the DCT and data-dependent matrices consistently improve as $k$ and $p$ increase. Quantitatively, in~\Cref{fig:traffic_all} and~\Cref{fig:shuttle_all}, we examine the relative error across all possible truncation and projection parameters. Through the difference color scales, we see that the identity and random transformations have larger relative errors overall and are more sensitive to the projection dimension than the truncation parameter (i.e., the error reduces predominantly along the $p$-axis). In contrast, the relative error for the DCT and data-dependent matrices is overall smaller and decreases when either $p$ or $k$ is increased.  There is a slightly slower decay of the relative errors for the \texttt{traffic} video (\Cref{fig:traffic_all}) than the \texttt{shuttle} video (\Cref{fig:shuttle_all}). This is because there is more high-frequency foreground activity in the \texttt{traffic} video (cars), and hence each successive increase in truncation and projection parameter makes incremental improvements. 

We capture this frequency information by examining the frontal slices in the transform domain in~\Cref{fig:traffic_features} and~\Cref{fig:shuttle_features}. For the \texttt{traffic} video, the transform domain features capture the movement of the vehicles at various frequencies. The transformations that perform best, the DCT and data-dependent matrices, separate the background from the foreground behavior effectively. In comparison, for the \texttt{shuttle} video, the transform domain features resemble the original spatial features for all choices of $\bfQ$. The features of higher-indexed frontal slices for the DCT and data-dependent matrices are negligible (almost zero). This is because the \texttt{shuttle} video has one main foreground object (rocket) moving at a fairly constant speed. As a result, the rocket's movement can be well-approximated by few frontal slices at appropriate frequencies. The identity and random matrices do not exploit this multilinear frequency information, and hence exhibit redundant features in the transform domain. 

In~\Cref{fig:traffic_matrix} and~\Cref{fig:shuttle_matrix}, we compare the top-performing transformation, $\bfU_3$, to the matrix SVD for various truncation parameters. In the left plot, we observe that for small projection dimensions $p$, the truncated matrix SVD produces a smaller relative error, but for larger $p$, the $\starQ$-SVD yields better approximations for all truncation parameters. For the full $p=n_3$ case, the $\starQ$-SVD provably gives a smaller error than the matrix SVD for the same truncation value~\cite[Theorem 5.3]{kilmer2019tensortensor}. In the right plot, we observe that the $\starQ$-SVD and matrix SVD are comparable in terms of compression ratio for small relative errors, and the $\starQ$-SVD offers slightly more compression for larger relative errors.  

Overall, the $\starQ$-representations with a good choice of transformation matrix efficiently capture multilinear behavior in both videos and are competitive with a matrix SVD representation with decreased storage costs for the less accurate approximations.

%% file: 05_01_video_figure.tex
\begin{figure}
\centering

\subfloat[Relative error $\|\TA - \TA_k(\bfQ)\|_F / \|\TA\|_F$  for truncation $k=5$ versus projection dimension $p$ for various choices of $\bfQ$. \label{fig:\myname_k5}]{
\begin{tikzpicture}
\draw[opacity=0] (-0.225\linewidth,0) -- (0.225\linewidth,0);
\node (n0) {\includegraphics[width=0.4\textwidth]{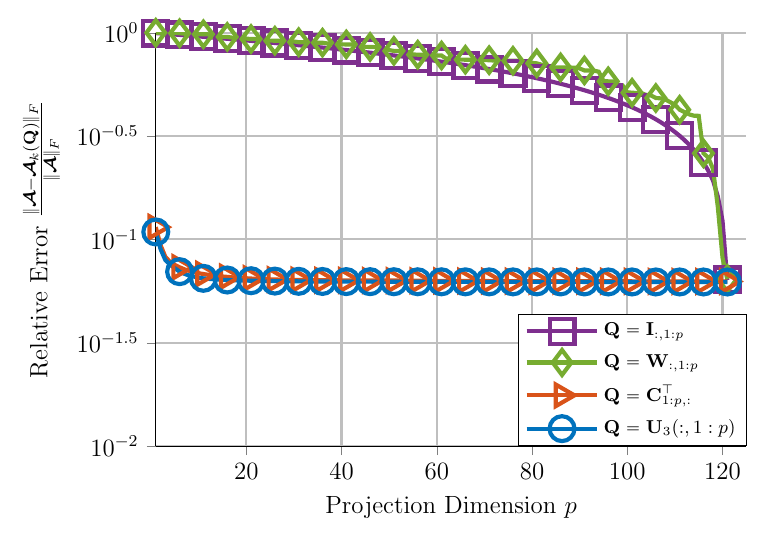}};
\end{tikzpicture}}
\hspace{0.01\linewidth}
\subfloat[Approximation of frame $45$ using $\TA_k(:,:,45)$ for each transformation. 
Each column corresponds to a different matrix and each row is a different combination of $k$ and $p$, increasing from top to bottom. \label{fig:\myname_approx}]{
\begin{tikzpicture}
\draw[opacity=0] (-0.25\linewidth,0) -- (0.25\linewidth,0);
\node (n0) {\includegraphics[width=0.5\textwidth]{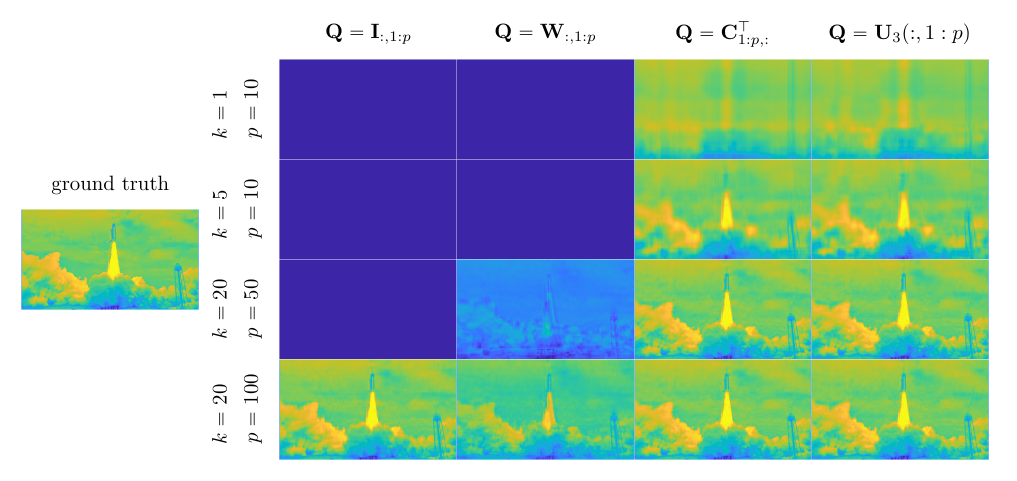}};
\end{tikzpicture}}

\subfloat[Relative error for all projection dimensions $p$ (rows) and truncations $k$ (columns) for all transformation matrices.  Colors close to black indicate smaller relative errors. Color scales are not equal in order to highlight patterns in approximation quality. \label{fig:\myname_all}]{
\begin{tikzpicture}
\draw[opacity=0] (-0.5\linewidth,0) -- (0.5\linewidth,0);
\node (n0) {\includegraphics[width=0.9\textwidth]{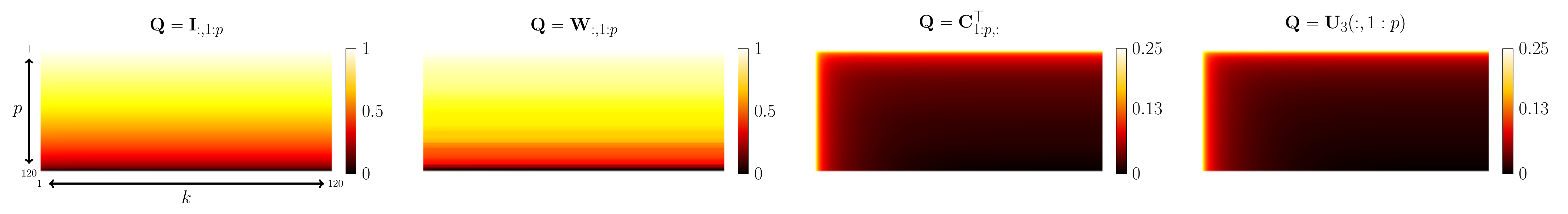}};
\end{tikzpicture}}

\subfloat[Features in the transform domain, where the slice $i$ is $\widehat{\TA}_{:,:,i} = \TA \times_3 \bfQ_{:,i}^H$.  
Each row corresponds to a different transformation using $p=n_3$ (no projection). 
Images are displayed at different color scales to highlight structure. \label{fig:\myname_features}]{
\begin{tikzpicture}
\draw[opacity=0] (-0.5\linewidth,0) -- (0.5\linewidth,0);
\node (n0) {\includegraphics[height=0.25\textwidth]{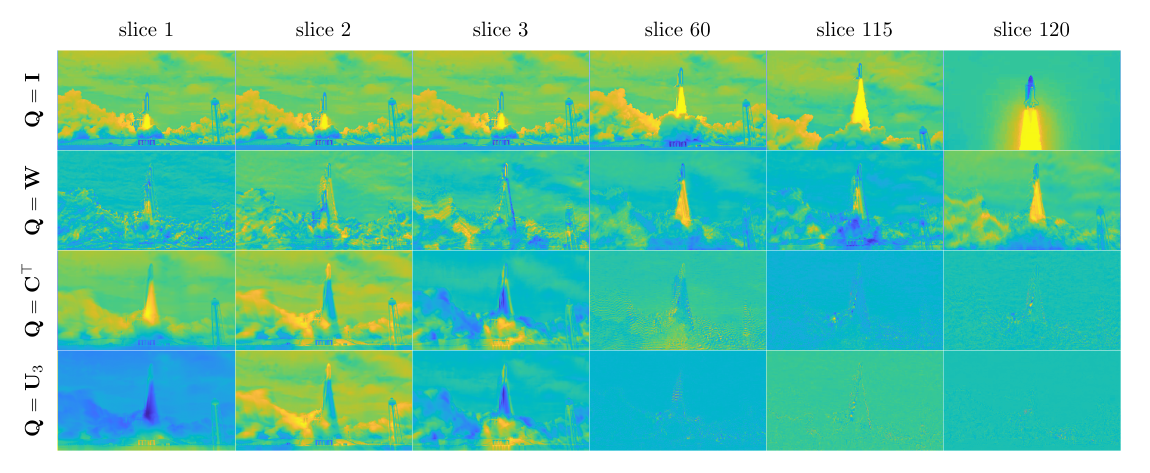}};
\end{tikzpicture}}

\subfloat[Relative error and compression ratio using $\bfQ = \bfU_3(:,1:p)$ for $p=10, 50, 100, n_3$ compared to the matrix SVD. (Left) Truncation versus relative error (lower is better) and (Right) relative error versus compression ratio (upper left is best). \label{fig:\myname_matrix}]{
\begin{tikzpicture}
\node (n0) {\includegraphics[width=0.4\textwidth]{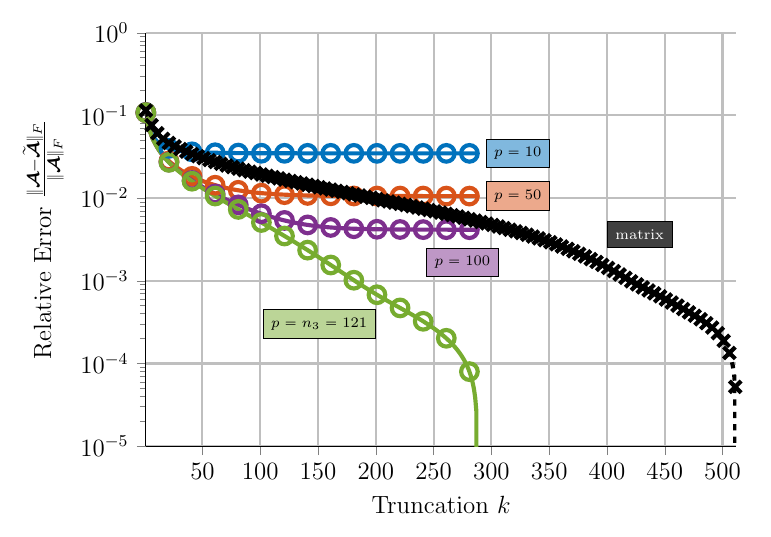}};
\node[above=0.0cm of n0.north, anchor=south, fill=lightgray, rounded corners=0.1cm, draw=black] (t0) {\scriptsize $k$ vs. RE};
\node[right=0.0cm of n0.east, anchor=west] (n1) {\includegraphics[width=0.4\textwidth]{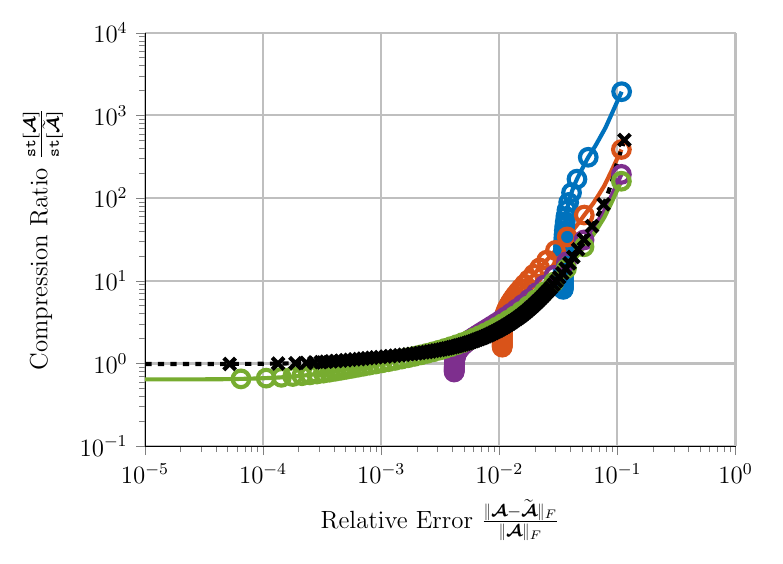}};
\node[fill=lightgray, rounded corners=0.1cm, draw=black] at (n1 |- t0) {\scriptsize  RE vs. CR};
\end{tikzpicture}}

\caption{Empirical results for \texttt{\myname} video compression.}
\label{fig:\myname}
\end{figure}

%% file: 05_02_hyperspectral.tex
\subsection{Hyperspectral Image Compression}
\label{sec:hyperspectral}
Hyperspectral images are naturally multilinear where each frontal slice corresponds to an image captured at a particular spectral bandwidth. 
In our experiments, we use the common Indian Pines dataset~\cite{IndianPines}, a $145 \times 145 \times 220$ tensor  of dimensions $\text{height} \times \text{width} \times \text{wavelength}$, available in the {\sc Matlab} Hyperspectral Toolbox\footnote{We can load the Indian Pines dataset using \texttt{hcube = hypercube('indian\_pines.dat'); A = hcube.DataCube;}} (see~\Cref{sec:hyperspectral_data} for a visualization). 
The two-dimensional renderings of the hyperspectral images in~\Cref{fig:hyperspectral1} are constructed by selecting three wavelengths as the RGB bands; in our case, we use (R,G,B) = (26, 16, 8).  
We examine the relative error and compression ratios for the $\starQ$-SVD and matrix SVD. 


\input{05_02_hyperspectral_figure}

\begin{figure}
    \centering
     \subfloat[Frobenius norm of each frontal slice of $\widehat{\TA} = \TA \times_3 \bfQ^H$ for each transformation using $p = n_3$. \label{fig:hyperspectral_frobenius}]{\includegraphics[width=0.4\linewidth]{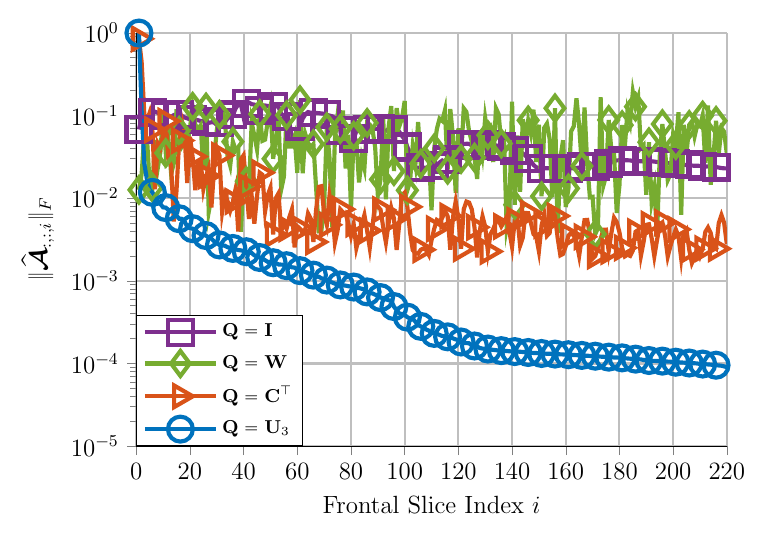}}
     \hspace{0.05\linewidth}
   \subfloat[Singular values of the mode-$3$ unfolding $\bfA_{(3)} = \bfU_3 \bfSigma_3 \bfV_3^H$ (\blue{blue}) and corresponding cumulative energy of the mode-$3$ unfolded tensor ({\color{red} red}). 
   	\label{fig:hyperspectral_decay}]{\includegraphics[width=0.4\linewidth]{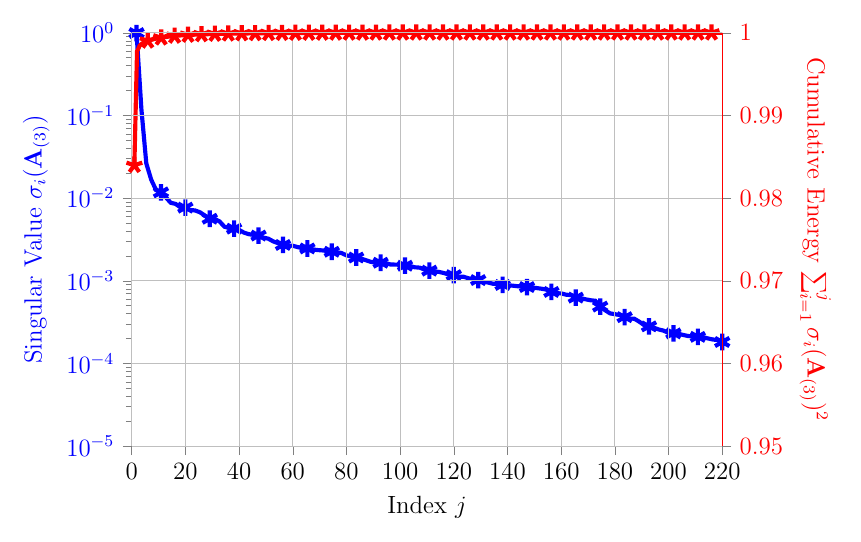}}

    \caption{Comparison of relative magnitudes of frontal slices in the transform domain (left) and mode-$3$ singular value energy (right). 
    The decay of the Frobenius norm for $\bfQ = \bfU_3$ exactly matches the decay of the mode-$3$ singular values.}
    \label{fig:hyperspectral_energy}
\end{figure}

In~\Cref{fig:hyperspectral_decay}, we see consistent evidence that the data-dependent matrix produces the best  $\starQ$-SVD approximation for the Indian Pines dataset. 
In~\Cref{fig:hyperspectral_rel_err}, we observe that the $\starQ$-SVD using $\bfQ = \bfU_3(:,1:p)$ achieves near optimal performance with about $1\%$ of the frontal slices to ($p\approx 2$ and $p / n_3 \approx 0.01)$. 
The DCT matrix performs second best, reaching near optimal performance using $p\approx 140$ ($p / n_3 \approx 0.64$). 
The RE to CR comparison plot shows that a relative error on the order of $10^{-1}$, the $\bfU_3$ representation requires almost three orders of magnitude less storage than the original data. 
In~\Cref{fig:hyperspectral_all}, we show the pattern of approximation quality for all combinations of $k$ and $p$.  
The data-dependent case obtains the smallest overall relative errors and its approximation quality improves similarly when increasing $p$ or $k$, with slightly more sensitivity to the choice of $k$. 
In comparison, the performance of the other three matrices is most sensitive to the choice of $p$, with the DCT matrix achieving significantly better approximations than the identity or random matrices. 
In~\Cref{fig:hyperspectral_approx}, we display approximations for different combinations of $k$ and $p$. 
Consistent with the relative error analysis, we observe that the $\bfU_3$ achieves the qualitatively accurate approximations even for low values of truncation and projection dimension. 
We note that when the identity and random matrices have poor approximations, the magnitude imbalance of the three color channels lead to unrealistic visualizations. 

In~\Cref{fig:hyperspectral_matrix}, we compare the $\starQ$-SVD with $\bfQ = \bfU_3(:,1:p)$ for different choices of $p$ to the matrix SVD. 
We see for small enough values of $p$, the $\starQ$-SVD achieves a better relative error to compression ratio performance. 
This provides numerical support about the advantages of leveraging multilinear correlations to form compressed representations. 

In~\Cref{fig:hyperspectral_energy}, we provide insight into why $\starQ$-SVD with $\bfQ = \bfU_3(:,1:p)$ admits the strongest performance for hyperspectral data compression. 
We examine the energy of the frontal slices of the tensor in the transform domain. 
We see in~\Cref{fig:hyperspectral_frobenius} that the norm of each transformed frontal slice $\|\widehat{\TA}_{:,:,i}\|_F$ for $i=1,\dots,n_3$ is at approximately the same magnitude for the identity and random transformations. 
This indicates that a quality $\starQ$-SVD representation has to approximate every frontal slice well, thereby requiring sufficiently large $k$ and $p$. 
In comparison, the DCT matrix shows more decay in the frontal slice magnitude, enabling better overall performance, even if the back frontal slices are poorly approximated. 
The data-dependent matrix shows the significant, monotonic decay in frontal slice magnitude which exactly matches the decay of the singular values of the mode-$3$ unfolding in~\Cref{fig:hyperspectral_decay} (see~\Cref{fig:mode3svd} for an explanation). 
Moreover, we see that almost all of the cumulative energy is captured in the first two mode-$3$ singular values. 
Thus, for $p \approx 2$, the $\starQ$-SVD using $\bfQ = \bfU_3(:,1:p)$ retains over $98\%$ of the overall energy.

Overall, the $\starQ$-SVD approximation using the data-dependent transformation is able to efficiently exploit high correlation among the tubes and can outperform the matrix SVD as a result.

%% file: 05_02_hyperspectral_figure.tex
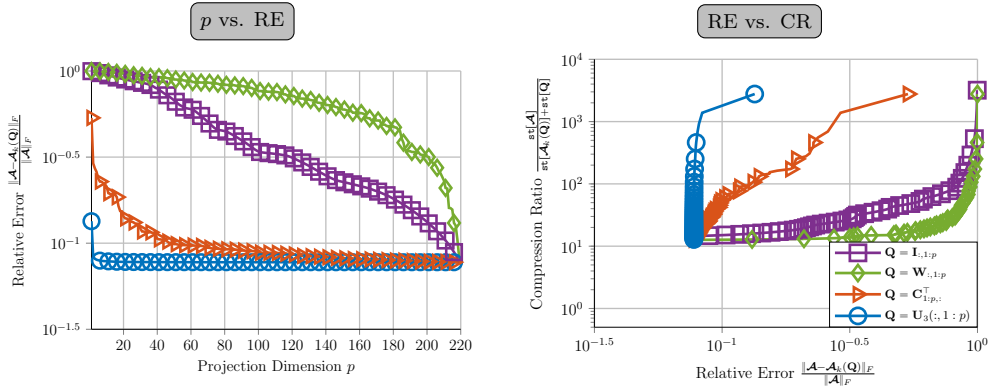
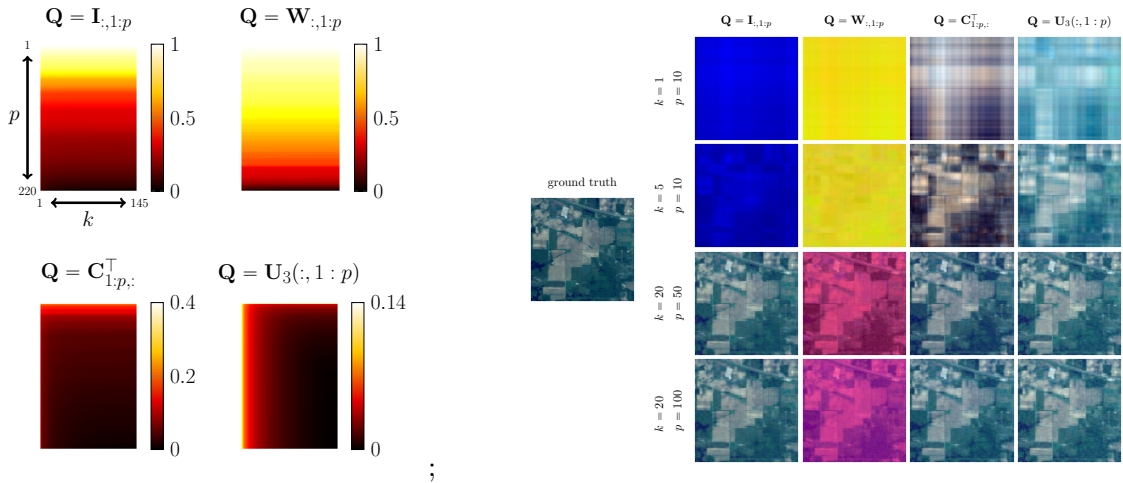
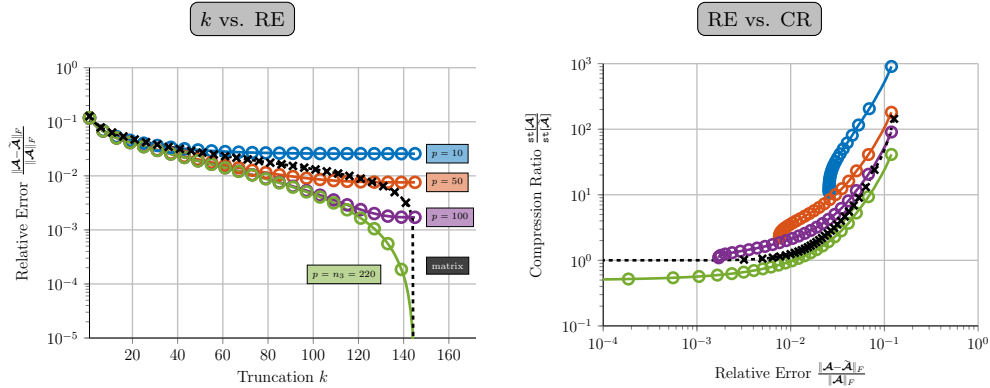
\begin{figure}
\centering

    \subfloat[Performance of the truncated $\starQ$-SVD with $k=5$ for all considered transformation matrices.  (Left) Projection dimension $p$ versus relative error (lower is better) and
    (Right) relative error versus compression ratio (upper left is better). \label{fig:hyperspectral_rel_err}]{\begin{tikzpicture}
       \node (n0) {\includegraphics[width=0.4\linewidth]{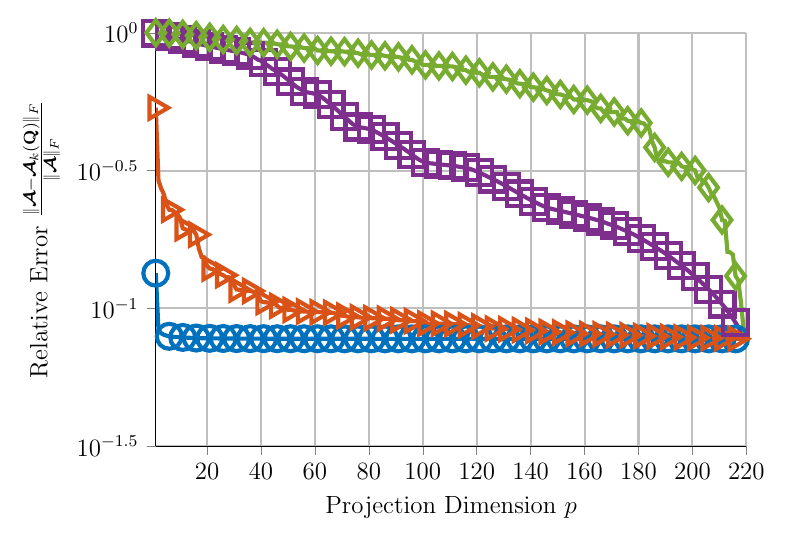}};
        \node[above=0.0cm of n0.north, anchor=south, fill=lightgray, rounded corners=0.1cm, draw=black] (t0) {\scriptsize $p$ vs. RE};
    	\node[right=0.0cm of n0.east, anchor=west] (n1) {\includegraphics[width=0.4\linewidth]{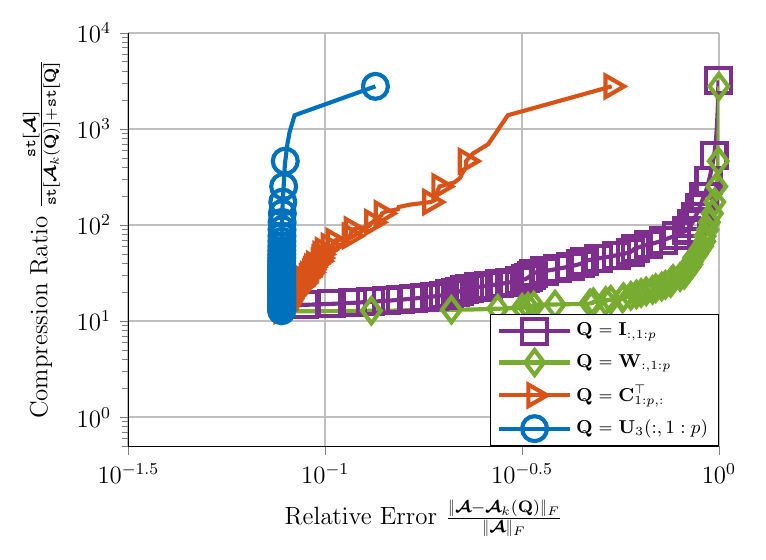}};
	\node[fill=lightgray, rounded corners=0.1cm, draw=black] at (n1 |- t0) {\scriptsize  RE vs. CR};
	\end{tikzpicture}
    }
   
    \subfloat[Relative error of the $\starQ$-SVD for all possible combinations of $p$ (rows) and $k$ (columns). 
    Colors closer to black indicate better performance. 
    The color scales are different per image to highlight structure. 
    \label{fig:hyperspectral_all}]{\includegraphics[width=0.35\linewidth]{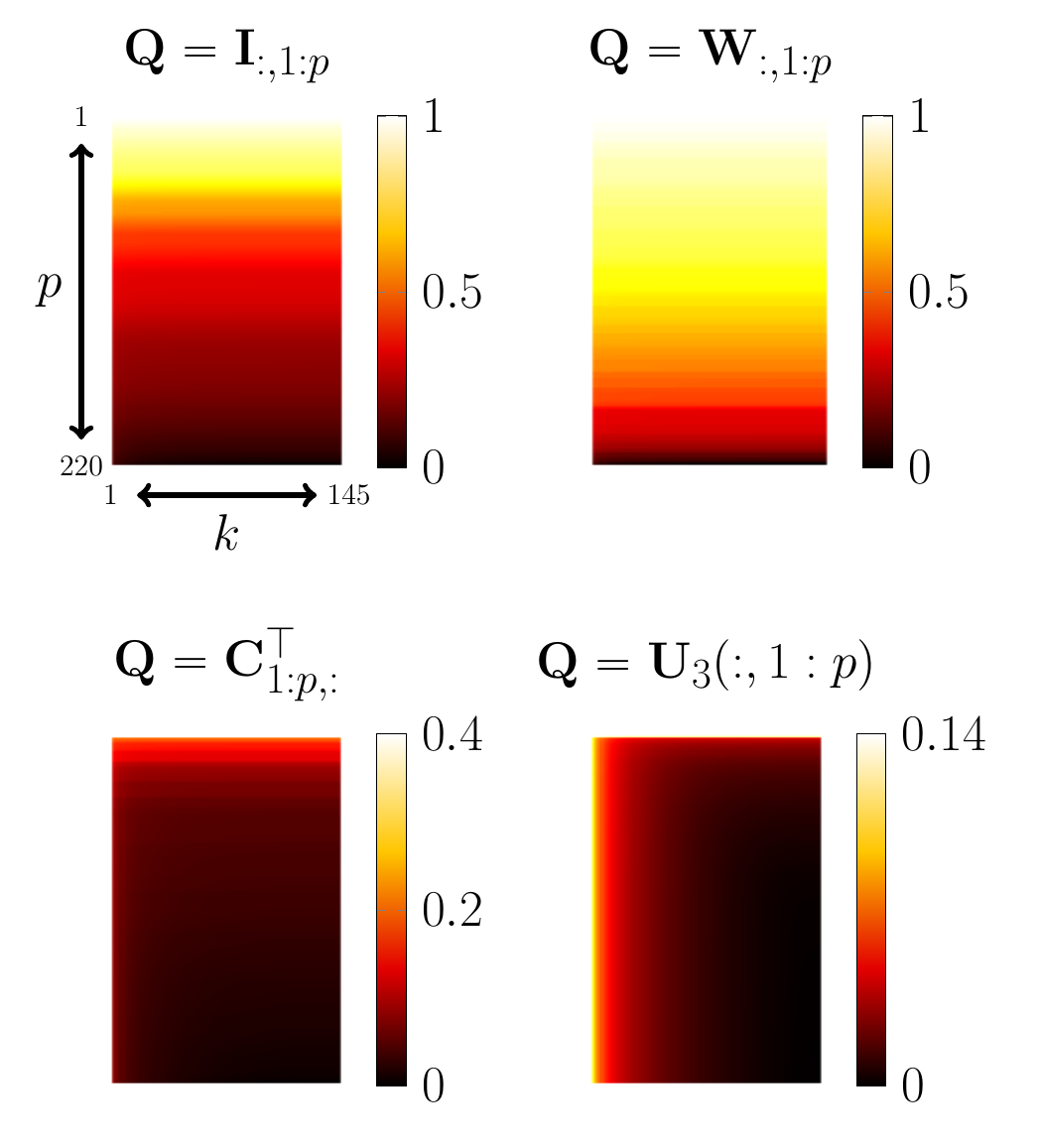}};
    \hspace{0.05\linewidth}
   \subfloat[Visualizations of the $\starQ$-SVD approximations to the hyperspectral data for various choices of $k$ and $p$. 
   Approximations are displayed as RGB images with channels (R,G,B) = (26, 16, 8) with rescaled pixel intensities between $0$ and $1$. 
   The bright blue, yellow, and magenta approximations indicate that certain color channels had significantly higher intensities than others. \label{fig:hyperspectral_approx}]{\includegraphics[width=0.5\linewidth]{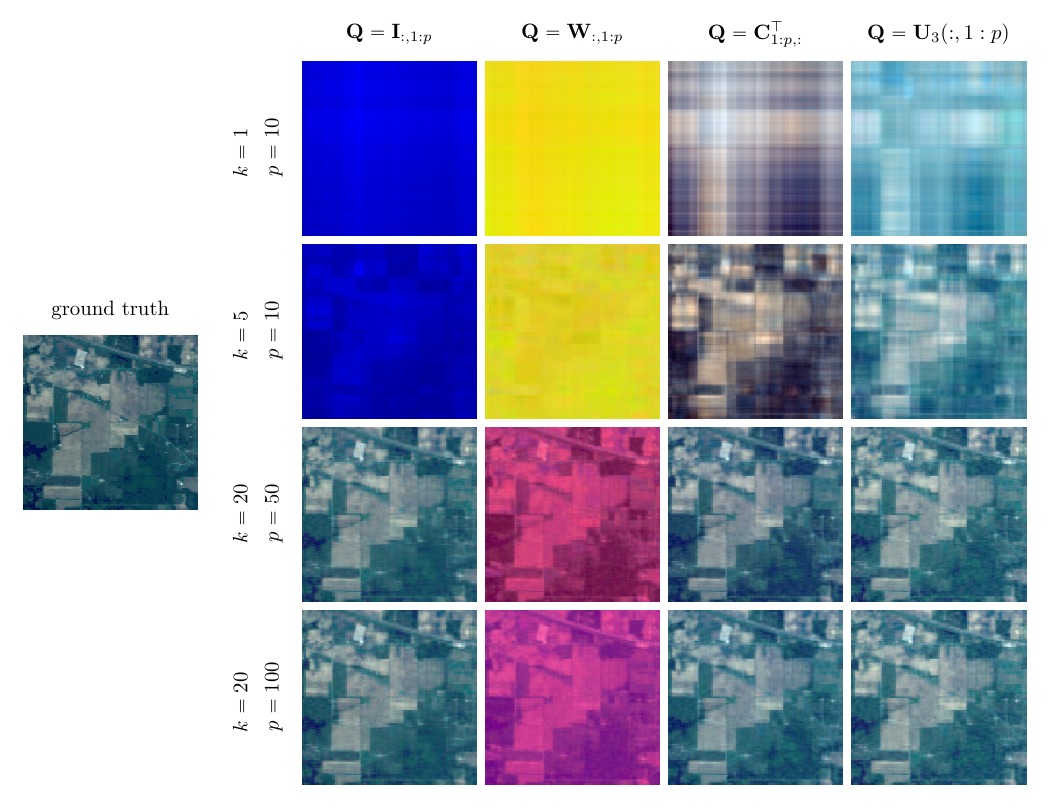}}

       \subfloat[Relative error and compression ratio using $\bfQ = \bfU_3(:,1:p)$ for $p=10, 50, 100, n_3$ compared to the matrix SVD. (Left) Truncation versus relative error (lower is better) and (Right) relative error versus compression ratio (upper left is best). \label{fig:hyperspectral_matrix}]{\begin{tikzpicture}
       \node (n0) {\includegraphics[width=0.4\linewidth]{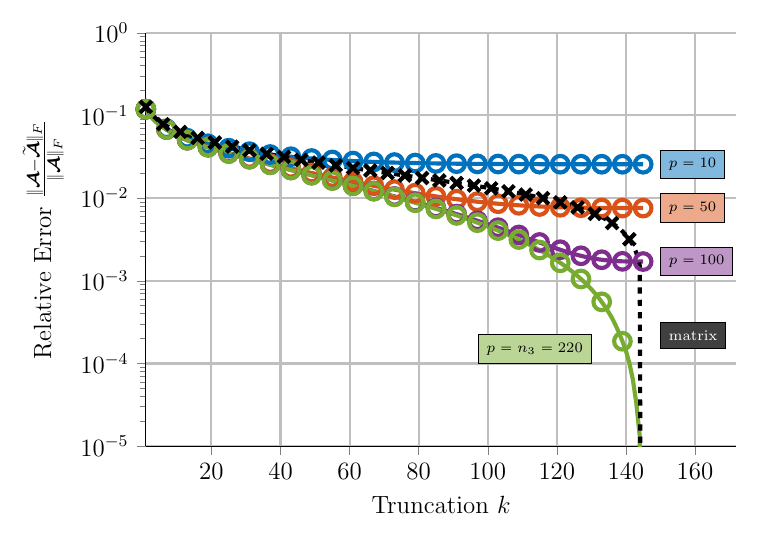}};
       \node[above=0.0cm of n0.north, anchor=south, fill=lightgray, rounded corners=0.1cm, draw=black] (t0) {\scriptsize $k$ vs. RE};
    	\node[right=0.0cm of n0.east, anchor=west] {\includegraphics[width=0.4\linewidth]{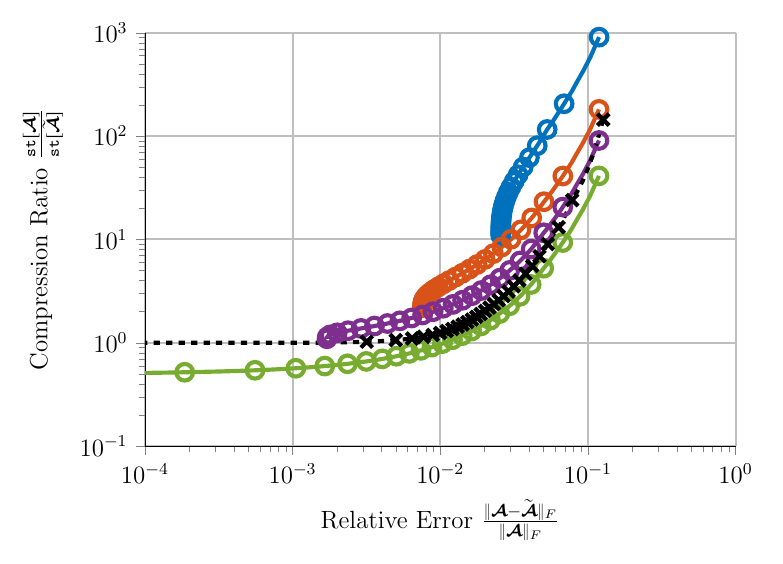}};
	\node[fill=lightgray, rounded corners=0.1cm, draw=black] at (n1 |- t0) {\scriptsize  RE vs. CR};
	\end{tikzpicture}
    }

\caption{Results for Indian Pines hyperspectral data compression.}
\label{fig:hyperspectral1}
\end{figure}

%% file: 05_03_hosvd.tex
\subsection{$\starQ$-SVD vs. HOSVD for Hyperspectral Image Compression}
\label{sec:tsvdq_vs_hosvd} 

\begin{figure}
\centering

\begin{tikzpicture}
\node (n0)  {\includegraphics[width=0.8\linewidth]{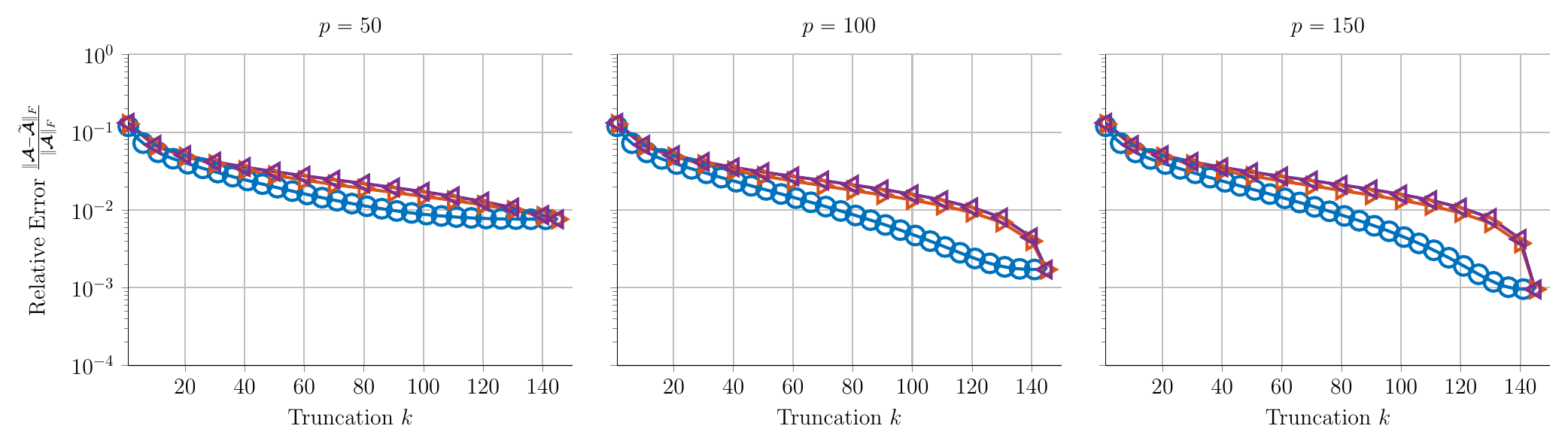}};
\node[below=0.0cm of n0.south, anchor=north] (n1) {\includegraphics[width=0.8\linewidth]{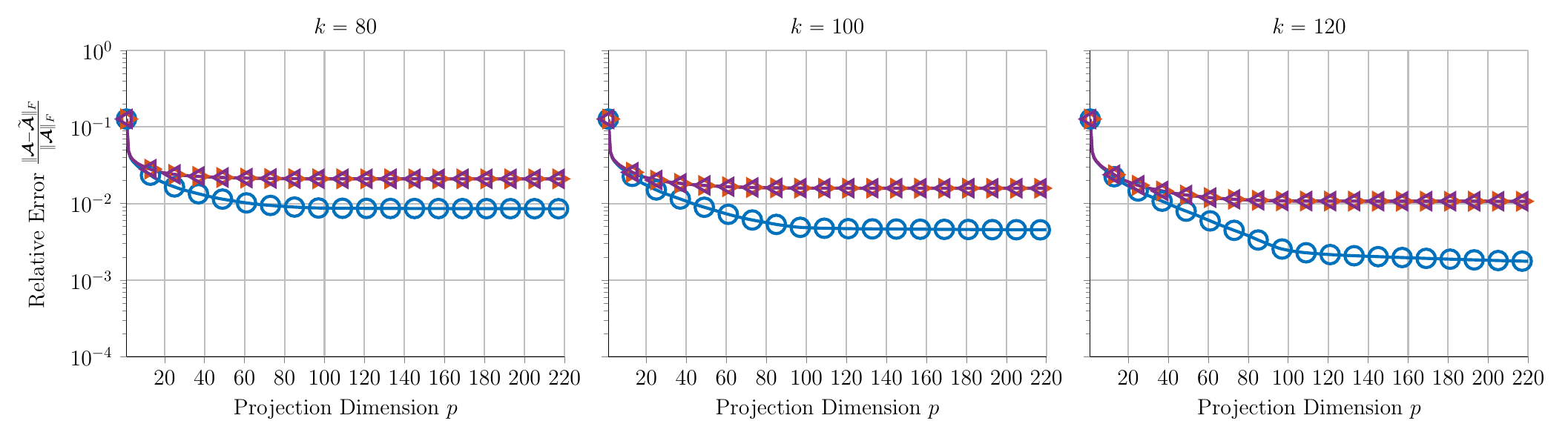}};
\node[below=0.0cm of n1.south, anchor=north] (n2) {\includegraphics[width=0.8\linewidth]{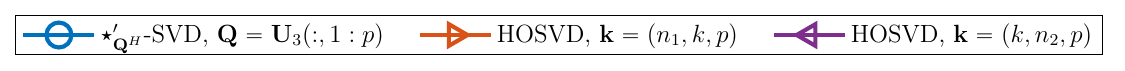}};
\end{tikzpicture}

\caption{Comparison of $\starQ$-SVD and HOSVD error for various choices of truncation $k$ and projection $p$ sizes. The top row varies the truncation parameter $k$ for fixed choices of $p=50,100,150$. The bottom row varies the projection dimension $p$ for fixed choices of $k=80,100,120$. To maximize HOSVD approximation performance, we choose the un-truncated dimension to be as large as possible. The $\starQ$-SVD obtains a smaller relative error for all choices of parameters, consistent with~\Cref{thm:projhosvd}.}
\label{fig:tsvdq_hosvd_err}
\end{figure}

In this section, we directly compare the $\starQ$-SVD, $\starQ$-SVDII, and HOSVD with $\bfQ = \bfU_3(:,1:p)$ for the Indian Pines hyperspectral dataset~\cite{IndianPines}. In~\Cref{fig:tsvdq_hosvd_err}, we empirically support~\Cref{thm:projhosvd} by comparing the relative error of the approximations for various choices of truncation and projection parameters $k$ and $p$, respectively. The $\starQ$-SVD achieves lower relative errors than the HOSVD for every combination of parameters $k$ and $p$, which provides empirical support for the theoretical bound in~\Cref{thm:projhosvd}. 

We now compare the more compressible the $\starQ$-SVDII to the truncated HOSVD. To ensure fair comparisons, we choose truncation parameters for each approximation starting from the $\starQ$-SVD truncation $k$. For the $\starQ$-SVDII (\Cref{alg:tsvdqII_multirank}), we choose the energy parameter $\gamma = \|\TA_k \times_3 \bfQ^H\|_F^2 / \|\TA \times_3 \bfQ^H\|_F^2$. This ensures that the $\starQ$-SVDII approximation will be no worse than the $\starQ$-SVD approximation (\Cref{thm:tsvdqII_vs_tsvdq}). For the truncated HOSVD, we consider two choices of multirank, $\bfk = (n_1,k_2,p)$ to achieve the best approximation or multirank $\bfk = (k_2,k_2,p)$ to achieve the most compression. We choose the HOSVD truncation parameter $k_2$ such that, when possible, the $\starQ$-SVD and HOSVD representations have similar storage costs.   Specifically, we approximately solve the following systems for the HOSVD truncation parameter $k_2$
\begin{subequations}\label{eq:hosvd_truncation}
\begin{align}
\text{HOSVD$(n_1,k_2,p)$:} & \quad k (n_1+n_2)p + n_3 p = n_1k_2p + n_1n_1 + n_2k_2 + n_3p\\
\text{HOSVD$(k_2,k_2,p)$:} & \quad k (n_1+n_2)p + n_3 p = k_2^2p + n_1k_2 + n_2k_2 + n_3p.
\end{align}
\end{subequations}
and round down to the closest positive integer value. The left-hand sides are the $\starQ$-SVD storage costs for a given $k$ and the right-hand sides are the HOSVD storage costs.  As a result of the choice of $k_2$, the truncated HOSVD will be no more expensive to store than as the $\starQ$-SVD.  

\begin{figure}
  \centering
 
  \includegraphics[width=0.8\linewidth]{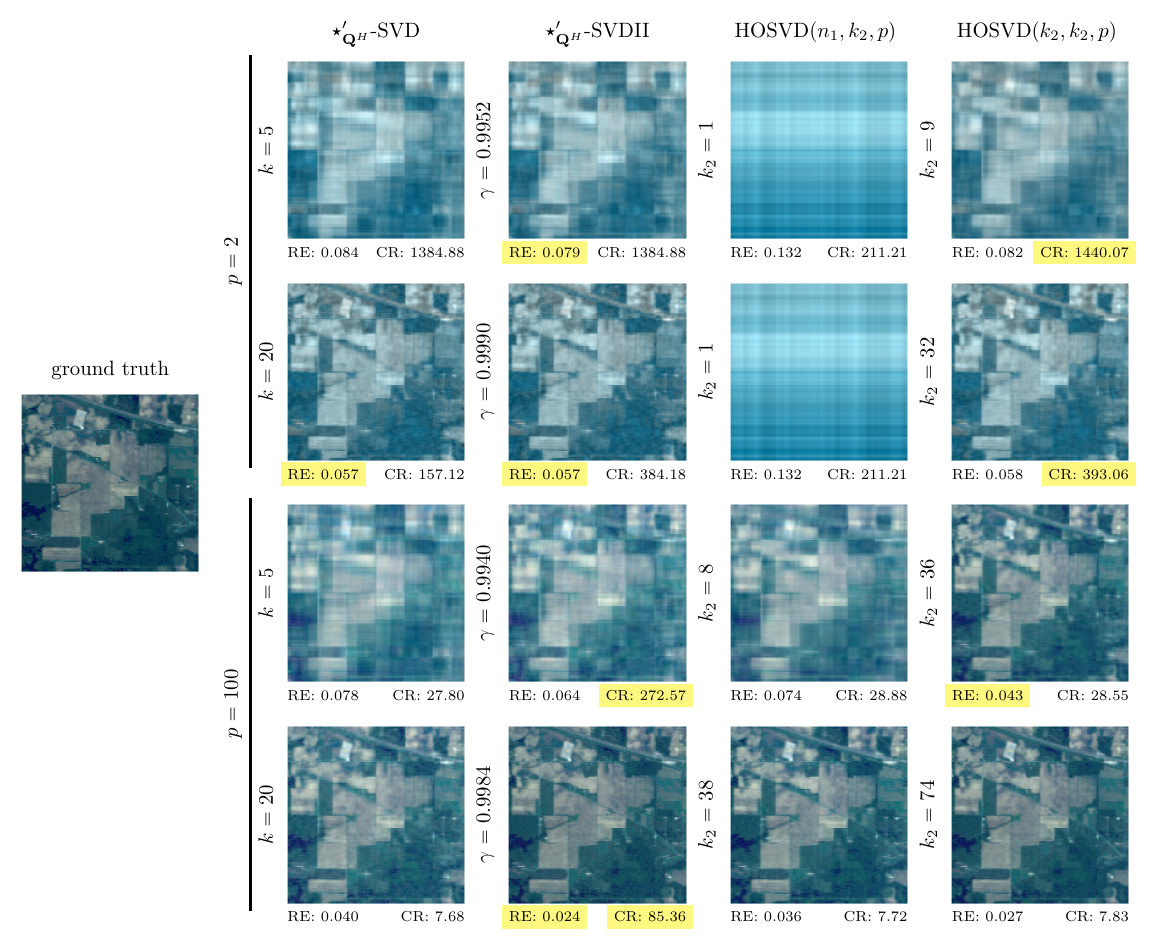}

 \caption{Approximations to the Indian Pines data using four different representations.  We fixed $p$ and $k$ for the $\starQ$-SVD, compute $\gamma$ for the $\starQ$-SVDII based on~\Cref{thm:tsvdqII_vs_tsvdq}, and choose $k_2$ such that the HOSVD examples have approximately the same amount of compression as the $\starQ$-SVD (if possible).  The choices of $p=2$  (top two rows) and $p=100$ (bottom two rows) roughly relate to the changes in singular value decay of the mode-$3$ unfolding depicted in~\Cref{fig:hyperspectral_decay}. The yellow color indicates the best relative error (RE) and compression ratio (CR) per row. }
  \label{fig:hyperspectral2_approx}
\end{figure}

We show the approximations for specific choices of truncations in~\Cref{fig:hyperspectral2_approx}. By construction, the $\starQ$-SVDII always always outperforms the $\starQ$-SVD with a lower relative error and higher compression ratio, as expected from~\Cref{thm:tsvdqII_vs_tsvdq}. Furthermore, the truncated HOSVD always has a larger compression ratio then the $\starQ$-SVD by design, except when there is no feasible choice of $k_2$ reach a higher level of compression. 
We observe empirically that the $\starQ$-SVDII and HOSVD approaches are competitive for the considered cases.

\begin{figure}
\centering

\subfloat[Fix $p$, vary $k$]{
\begin{tikzpicture}
\node  (n0) {\includegraphics[width=0.48\linewidth]{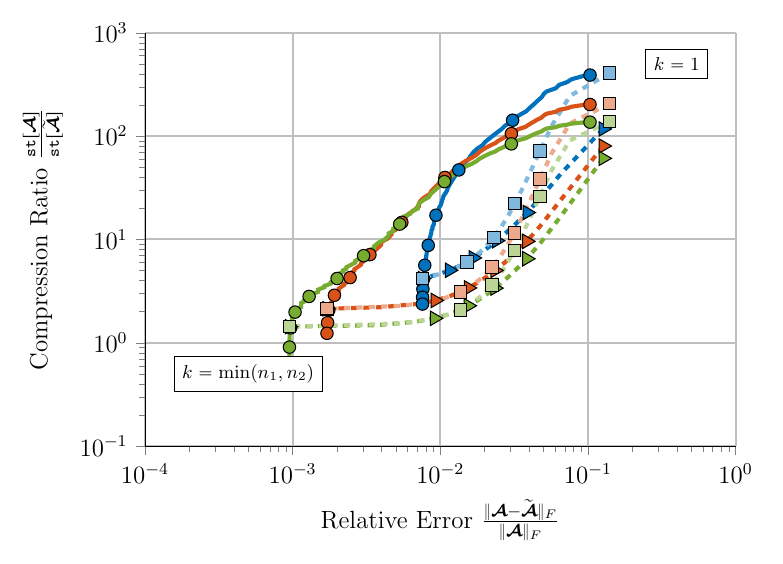}};
\node[above=-0.5cm of n0.north, anchor=south] {\includegraphics[width=0.38\linewidth]{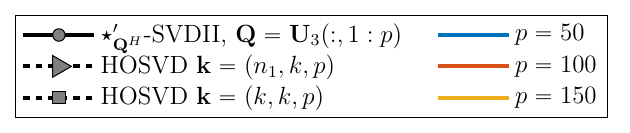}};
\end{tikzpicture}
}
\subfloat[Fix $k$, vary $p$]{
\begin{tikzpicture}
\node  (n0) {\includegraphics[width=0.48\linewidth]{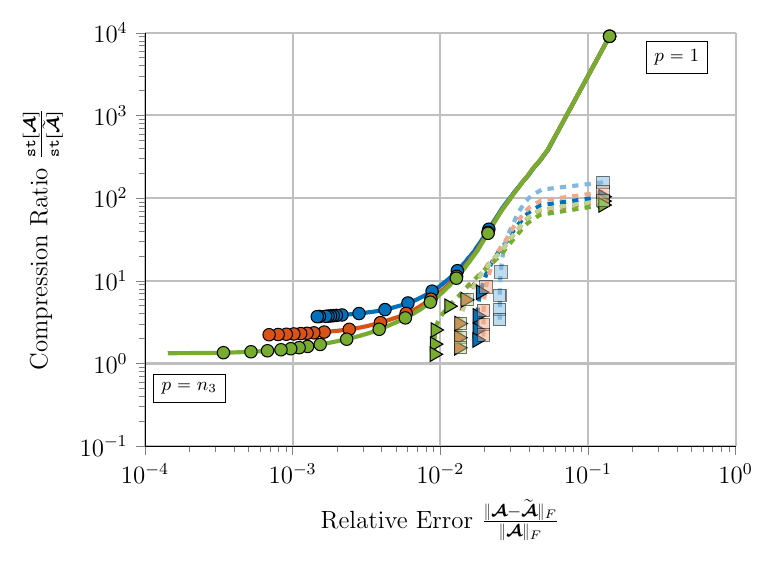}};
\node[above=-0.5cm of n0.north, anchor=south] {\includegraphics[width=0.38\linewidth]{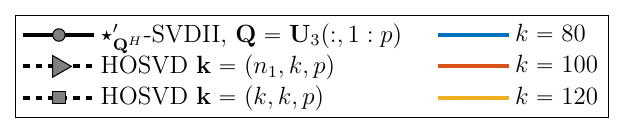}};
\end{tikzpicture}
}

\caption{Relative error versus compression ratio to compare the $\starQ$-SVDII and the truncated HOSVD with various multilinear rank combinations. We denote any approximation as $\widetilde{\TA}$. Compression strategies with lower relative error and a higher compression ratio are better (upper left is best). Solid lines with circle markers are used for the $\starQ$-SVDII, dashed lines with darker colors and triangle markers are used for the HOSVD$(n_1,k,p)$ cases, and dashed lines with lighter colors and square markers are used for the HOSVD$(k,k,p)$ cases. The two HOSVD cases represent the best possible relative error and best possible compression ratio, respectively, under the assumptions of~\Cref{thm:projhosvd}. The left plot varies the truncation parameter $k$ for three choices of $p = 50,100,150$ and the right plot varies the projection parameter $p$ for three choices of truncation $k=80,100,120$.}
\label{fig:tsvdqII_vs_hosvd} 
\end{figure}

To compare the $\starQ$-SVDII to the HOSVD performance across parameters, we plot the relative error versus compression ratio in~\Cref{fig:tsvdqII_vs_hosvd}. As before, we generate an energy parameter for the $\starQ$-SVDII based on the truncation $k$ for the $\starQ$-SVD. For the HOSVD, we vary the truncation $k =1,10,20,...,140, 145$ where $k = 145$ represents the non-truncated case. We see that for various choices of $k$ and $p$, the $\starQ$-SVDII achieves better approximation quality for less storage. We observe that both HOSVD truncation strategies have similar performance because the $\bfk = (n_1,k_2,p)$ case minimizes the relative error whereas the $\bfk = (k_2,k_2,p)$ case maximizes the compression ratio. 

The main takeaway is that the $\starQ$-SVDII consistently outperforms the truncated HOSVD in terms of relative error to compression ratio. This demonstrates that the projected product can have both theoretical and numerical advantages over the HOSVD, making it an appealing tensor representation strategy for data with highly-correlations along the third dimension.

%% file: 06_conclusions.tex
\section{Conclusions}
\label{sec:conclusions}

We developed a unified algebraic framework for projected tensor-tensor products that preserves matrix mimeticity with reduced computational overhead. In~\Cref{sec:projected_products}, we verified that fundamental linear algebraic properties are inherited under the projected product. In~\Cref{sec:tsvdq}, we introduced the $\starQ$-SVD, proved an Eckart-Young-like optimality result, and provided insight into an optimal choice of transformation matrix. We extended the $\starQ$-SVD to the more compressible $\starQ$-SVDII in~\Cref{sec:tsvdqII} and proved the $\starQ$-representational superiority over the truncated HOSVD in~\Cref{sec:hosvd}. In~\Cref{sec:numericalexperiments}, we provided extensive experiments on video and hyperspectral data to empirically verify the theoretical findings and demonstrate the compressibility of $\starQ$-representations. We observed that the $\starQ$-SVD (and variants) with an appropriate choice of transformation produces the best representation in terms of relative reconstruction error compared to compression ratio in all of our experiments.

Our work lays the foundation for future exploration of other tensor algebras defined by non-invertible, efficient transformations. 
In the short term, we can build from the work in~\cite{newman2024optimalmatrixmimetictensoralgebras} to optimize the transformation matrix $\bfQ$ practically and we can extend to higher-order projected products following the work in~\cite{keegan2022a}. 
As a subsequent extension, we can generalize to any linear transformation by using the Moore-Penrose pseudoinverse~\cite{Penrose_1955, Barata_2011} to approximate the reverse transformation. 
This will enable exploration into a wider range of matrix structures that can be stored and (pseudo) inverted efficiently, e.g., low-rank, symmetric positive semidefinite matrices, but may sacrifice some algebraic guarantees. 
Beyond algebraic extensions, we will explore new applications to compress dense tensors with at least one large dimension, such as spatio-temporal data for chemo-sensing~\cite{VERGARA2013462} and fluid flow simulations~\cite{cabot_reynolds_2006, zhao_sdrbench_2020}.

%% file: 07_projected_products_example.tex
\section{Projected Products Example}
\label{sec:projected_products_example}

We provide an example to illustrate the various perspectives of projected products in connection with the properties presented in~\Cref{sec:equivalent_presentations}. 
Consider the real-valued tubes $\bfa, \bfb\in \Rbb^{1\times 1\times 4}$ with
	\begin{align}
	\bfa = \boxed{\begin{array}{cccc} 2 & 4 & 6 & 8 \end{array}} \qquad \text{and} \qquad
	\bfb = \boxed{\begin{array}{cccc}  1 & -1 & 1 & 0\end{array}} .
	\end{align}
Let $\bfM = \bfH_4^\top$, the transposed of the $4\times 4$ Haar wavelet matrix given by
	\begin{align}
	\bfH_4 = \frac{1}{2}\begin{bmatrix}
	1 & 1 & 1 & 1\\
	1 & 1 & -1 & -1\\
	\sqrt{2} & -\sqrt{2} & 0 & 0\\
	0 & 0 & \sqrt{2} & -\sqrt{2}
	\end{bmatrix}.
	\end{align}
Partition $\bfM = \begin{bmatrix} \bfQ & \bfQ_{\perp}\end{bmatrix}^H$ where $\bfQ\in \Rbb^{4\times p}$ and $p = 2$. 
Then, 
	\begin{alignat}{3}
	&\bfa \starM \bfb && &&=  \boxed{\begin{array}{cccc}   2 & 4 & 3 & 8 \end{array}} \label{eq:abM}\\
	&\bfa \starQ \bfb &&= (\bfa \times_3 \bfQ \bfQ^H) \starM (\bfb \times_3 \bfQ \bfQ^H) 
	&&= \boxed{\begin{array}{cccc}   3 & 3 & 3 & 0 \end{array}} \label{eq:abQ1}\\
	&\bfa \starQperp \bfb 
	&&= (\bfa \times_3 \bfQ_{\perp} \bfQ_{\perp}^H) \starM (\bfb \times_3 \bfQ_{\perp} \bfQ_{\perp}^H) &&= \boxed{\begin{array}{cccc}  -1 & 1 & 0 & 8\end{array}} \label{eq:abQ2}
	\end{alignat}
where $\bfQ_{\perp} \bfQ_{\perp}^H = \bfI_4 - \bfQ \bfQ^H$.  
The projected tubes are given by
	\begin{align}
	\bfa \times_3 \bfQ \bfQ^H &= \boxed{\begin{array}{cccc}   3 & 3 & 6 & 0  \end{array}}
	&
	\bfb \times_3 \bfQ \bfQ^H &=\boxed{\begin{array}{cccc}   0 & 0 & 1 & 0  \end{array}}\label{eq:Q}\\
	\bfa \times_3\bfQ_{\perp} \bfQ_{\perp}^H &=\boxed{\begin{array}{cccc}  -1 & 1 & 0 & 8  \end{array}}
	&
	\bfb \times_3 \bfQ_{\perp} \bfQ_{\perp}^H &=\boxed{\begin{array}{cccc}   1 & -1 & 0 & 0   \end{array}}. \label{eq:Qperp}
	\end{align}
From~\eqref{eq:abM},~\eqref{eq:abQ1}, and~\eqref{eq:abQ2}, we see that $\bfa \starM \bfb = \bfa \starQ \bfb + \bfa \starQperp \bfb$, as expected from~\eqref{eq:colSpaceDef}. 
We further see that the tubes in~\eqref{eq:Q} lie in the column space of $\bfQ$ and the null space of $\bfQ_{\perp}^H$, following the property in~\eqref{eq:nullspacePerspective}.  
Applying $\bfM$ along the projected tubes, we obtain the properties from~\eqref{eq:backFrontalSlicesZero}, though the numbers are not nice enough to be worth showing.

%% file: 07_numerical_experiment_data.tex
\section{Numerical Experiment Data}
\label{sec:numerical_experiments_data}

\subsection{Video Data}
\label{sec:video_data}

\begin{figure}
\centering

\subfloat[\texttt{traffic}]{\includegraphics[width=0.4\linewidth]{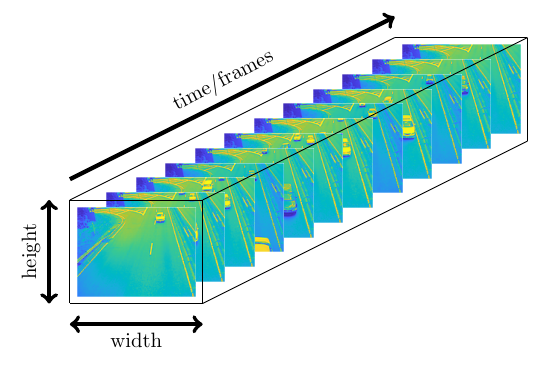}}
\subfloat[\texttt{shuttle}]{\includegraphics[width=0.4\linewidth]{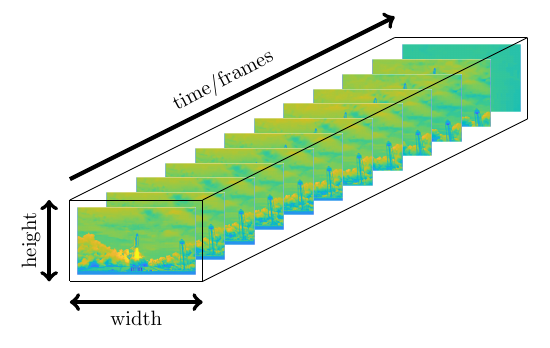}}

\caption{Orientation of video datasets. The images are not drawn to scale nor are the color scales equal.}
\label{fig:video_data}
\end{figure}

We present a visualization of the two video datasets, \texttt{traffic} and \texttt{shuttle}, in~\Cref{fig:video_data}. 
Our goal is to present the orientation of the data for transparency of the experiments. 

\subsection{Hyperspectral Data}
\label{sec:hyperspectral_data}

\begin{figure}
\centering

\subfloat[Hyperspectral cube with three spatial locations indicated by the different colored circles.]{\includegraphics[width=0.32\linewidth]{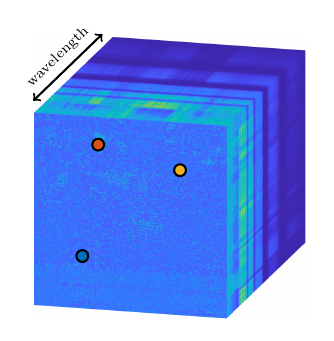}}
\hspace{0.05\linewidth}
\subfloat[Example of spectral patterns at three locations and a color visualization of the data using channels RGB = (26,16,8). The three channels are shown as vertical dashed lines. 
The normalized pixel intensities of the three locations are well-separated at these bands, enabling the various regions to be clearly delineated.]{\includegraphics[width=0.58\linewidth]{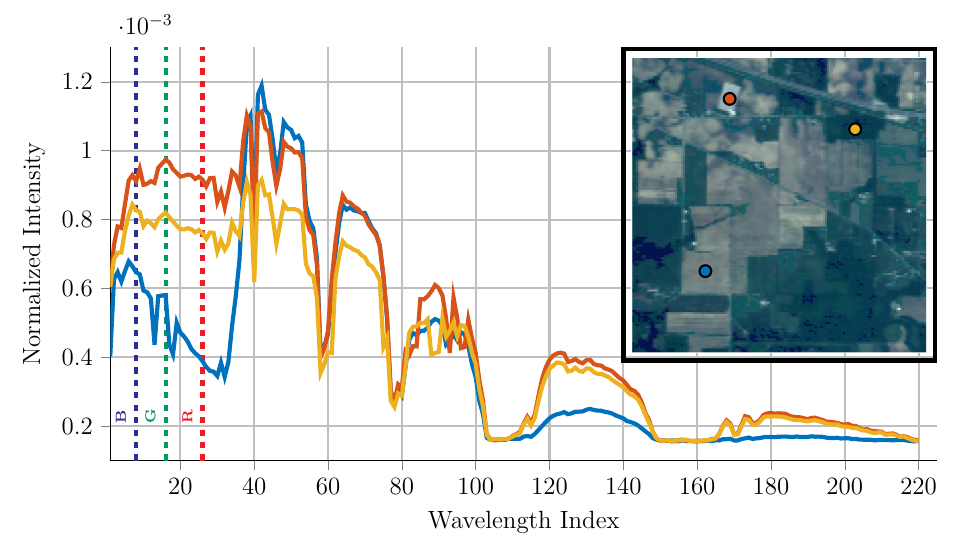}}

\caption{Visualization of Indian Pines dataset from the {\sc Matlab} Hyperspectral Imaging Toolbox. 
The colors of the dots in the left picture correspond to the three solid lines in the right picture. }
\label{fig:hyperspectral_data}
\end{figure}

We present a visualization of the Indian Pines dataset~\cite{IndianPines} in~\Cref{fig:hyperspectral_data}. 
The Indian Pines data depicts a birds-eye view of a region in North-western Indiana that consists of forests, crops, and man-made infrastructure. 
The various features of the region illicit different spectral signatures. 

